\documentclass[10pt]{article}

\usepackage{amsmath, amssymb, amsfonts, amsthm, graphicx, overpic}

\newcommand{\ptype}[2]{\genfrac{}{}{0pt}{}{#1}{#2}}
\newcommand{\MR}{\mathsf{MR}}
\newcommand{\PG}{\mathsf{PG}}
\newcommand{\Z}{\mathbb{Z}}
\newcommand{\D}{\mathbb{D}}
\newcommand{\Proj}{\mathbb{P}}
\newcommand{\collection}[1]{\mathsf{#1}}
\newcommand{\CS}{\collection{S}}

\theoremstyle{plain}
\newtheorem{theorem}{Theorem}
\newtheorem{lemma}{Lemma}
\newtheorem{corollary}{Corollary}
\newtheorem{proposition}{Proposition}
\newtheorem{problem}{Problem}
\newtheorem{claim}{Claim}
\newtheorem*{mrtheorem}{Motzkin-Rabin Theorem}
\newtheorem*{sgtheorem}{Sylvester-Gallai Theorem}

\theoremstyle{definition}
\newtheorem{example}{Example}

\begin{document}

\bibliographystyle{amsplain}

\title{Blocking sets in small finite linear spaces\thanks{This material is based upon work supported by the National Research Foundation under Grant number 2053752.}}
\author{Lou M. Pretorius\thanks{Department of Mathematics and Applied Mathematics,
        University of Pretoria,
        Pretoria 0002, South Africa.
        E-mail: \texttt{lpretor@scientia.up.ac.za}}
\and
        Konrad J. Swanepoel\thanks{Department of Mathematics, Applied Mathematics 
        and Astronomy, University of South Africa,
        PO Box 392, UNISA 0003, South Africa.
        E-mail: \texttt{swanekj@unisa.ac.za}}}
\date{}
\maketitle

\vspace{-7mm}
\begin{abstract}
We classify all finite linear spaces on at most $15$ points admitting a blocking set.
There are no such spaces on $11$ or fewer points, one on $12$ points, one on $13$ points, two on $14$ points, and five on $15$ points.
The proof makes extensive use of the notion of the weight of a point in a $2$-coloured finite linear space, as well as the distinction between minimal and non-minimal $2$-coloured finite linear spaces.
We then use this classification to draw some conclusions on two open problems on the $2$-colouring of configurations of points.
\end{abstract}

\section{Introduction}

\subsection{Finite linear spaces and blocking sets}\label{subsec:intro}
By a \emph{finite linear space} (or f.l.s.) $\CS$ we mean a triple consisting of
a finite set of \emph{points} (denoted by capital letters $A, B, \dots$), a finite set of \emph{lines} (denoted by lower case letters $\ell, m, \dots$), and an \emph{incidence relation} $\in$ between the set of points and the set of lines ($P\in\ell$ being read as ``$P$ \emph{lies on} $\ell$'' or ``$\ell$ \emph{passes through} $P$''), such that the following properties hold:
\begin{description}
\item[I1] For any two distinct points $P$ and $Q$ there is a unique line passing through $P$ and $Q$ (denoted by $PQ$).
\item[I2] Each line passes through at least two points.
\item[I3] There exist three non-collinear points.
\end{description}
By I1 and I2 we may identify a line with the set of points lying on it.
By I3 we do not allow a finite linear space to be ``one-dimensional''.
Of course we could have omitted I3 (as is done by certain authors) and then we would only have to change the formulations of certain statements.
Other equivalent terms used in the literature are simple matroid of rank $3$, and pairwise balanced incomplete block design.
A \emph{proper finite linear space} is an f.l.s.\ with each line passing through at least $3$ points.
Proper finite linear spaces are important when considering generalizations of the Sylvester-Gallai Theorem (see Section~\ref{subsec:sg}).
A proper f.l.s.\ is also called an SG design in \cite{MR47:3207}.
Our terminology follows recent papers such as \cite{MR95j:51018, MR2000f:05020, MR2002g:05042}.

A \emph{blocking set} of a finite linear space $\CS$ is a subset $B$ of $\CS$ such that each line of $\CS$ contains a point of $B$ and a point of $\CS\setminus B$.
There is extensive literature on blocking sets of block designs (see \cite{MR95j:05028, MR1178507} and their references), and in particular of affine and projective planes (see \cite[Chapter~8]{MR99c:51001}, \cite[Chapter~13]{Hirschfeld}, \cite{BE, Bl,Sz-Sz} and their references).
Not as much has been done on blocking sets for general linear spaces.
Cameron \cite{MR87e:51014} mentions in passing that blocking sets may be studied in general finite linear spaces.
Hahn has two papers \cite{MR92h:51016, MR94a:51017} on how the existence of certain types of blocking sets in finite linear spaces forces the space to be almost a projective or affine plane.
Hoffman, Lindner and Phelps \cite{MR92f:05016} and Franek, Griggs, Lindner and Rosa \cite{MR2003a:05026} show that a linear space on $v$ points with four points on each line that admits a blocking set exists whenever $v\equiv1, 4\pmod{12}$.
Batten, Coolsaet and Street \cite{MR97i:05010} consider the existence of blocking sets in finite linear spaces where each line passes through either $2$ or $4$ points.
They show that apart from the above-mentioned cases $v\equiv1, 4\pmod{12}$, such spaces also exist for $v\equiv 2\pmod{12}$, and do not exist for all remaining congruence classes, except perhaps $v\equiv 10\pmod{12}$, where they conjecture that no such finite linear spaces admit blocking sets.
Ling \cite{MR2000m:05028} proves that an f.l.s.\ with five points on each line that admits a blocking set exists whenever $v\equiv1,5\pmod{20}$.
See also \cite{MR93m:05148, MR98d:05023, MR2003f:05026} for related results.

\smallskip
Our main result is a classification of all finite linear spaces of at most $15$ points that admit blocking sets.
\begin{theorem}\label{mainthm}
\mbox{}\begin{enumerate}
\item No finite linear space on $\leq 11$ points admits a blocking set.
\item The punctured projective plane of order $3$ is the only f.l.s.\ on $12$ points with a blocking set.
\item The projective plane of order $3$ is the only f.l.s.\ on $13$ points with a blocking set.
\item There are exactly two finite linear spaces on $14$ points with blocking sets.
Both are subsets of the projective plane of order $4$.
\item There are exactly five finite linear spaces on $15$ points with blocking sets.
All except one are subsets of the projective plane of order $4$, with the exceptional one not embeddable in any Desarguesian projective space.
Furthermore, all except one are extensions of one of the two $14$-point spaces with blocking sets.
\end{enumerate}
\end{theorem}
It is rather surprising that there are so few finite linear spaces with blocking sets up to $15$ points.
In all of them the blocking sets are unique up to complements and automorphisms of the space, with the exception of one of the $15$-point spaces which admits two essentially different blocking sets.
In Section~\ref{examples} we describe the finite linear spaces mentioned in the above theorem and then give a precise formulation of the theorem (see Theorem~\ref{mainthm2}).
The proof of the theorem is an extensive case analysis, as is to be expected.
The two main tools used in the proof are the notions of \emph{minimal} two-coloured f.l.s.\ (see Section~\ref{subsec:def}) and of the \emph{weight} of a point in a two-coloured f.l.s.\ that we introduce below (see Lemma~\ref{weight} and the definition preceding it).

Note that finite linear spaces up to $12$ points have been enumerated \cite{MR99m:05034}, with $13$ points still out of reach.
A brute-force computer proof of our result is therefore not yet feasible.

\subsection{General results}
We now give an overview of general results on the existence of blocking sets which we use in the proof of our classification.

It is well known that the projective plane of order $3$ has a blocking set (see Example~\ref{ex:thirteen}).
It is easily seen that the punctured projective plane of order $3$ still has a blocking set (Example~\ref{ex:twelve}).
This is a minimal example: we show that an f.l.s.\ on at most $11$ points does not admit a blocking set.
More generally, if an f.l.s.\ admits a blocking set, the blocking set must have size $\geq 6$ (Lemma~\ref{atleastsixpoints}).
Also, for each $v\geq 12$, there exists an f.l.s.\ on $v$ points admitting a blocking set (Proposition~\ref{exist}).
An f.l.s.\ with a blocking set of $6$ points has at most $13$ points (Proposition~\ref{sixpoints}).

It is a folklore result that an f.l.s.\ with at most three points on a line does not admit a blocking set.
Thus a necessary condition for an f.l.s.\ to admit a blocking set is that the f.l.s.\ has a line with at least $4$ points.
A further necessary condition is that there are at least $4$ lines passing through \emph{each} point (Lemma~\ref{atleastfourlines}), and more generally, that the f.l.s.\ is not the union of $3$ lines.

\subsection{MR geometries: Definitions}\label{subsec:def}
Since the complement of a blocking set is also a blocking set, we find it more convenient to formulate our results in terms of $2$-colourings.
A \emph{$2$-colouring} $\chi$ of the finite linear space $\CS$ is a function assigning the colour red or green to each point of $\CS$.
In all the figures we use $\circ$ for red and $\bullet$ for green.
A line is \emph{monochromatic} if all points incident with it have the same colour.
A \emph{proper \textup{(}weak\textup{)}} $2$-colouring $\chi$ is a $2$-colouring with the property that no line is monochromatic.
An \emph{MR geometry} $(\CS,\chi)$ is a properly $2$-coloured finite linear space.
Thus each colour class of an MR geometry $(\CS,\chi)$ forms a blocking set of the underlying finite linear space $\CS$.
Conversely, any blocking set of $\CS$ determines a proper $2$-colouring.
Thus an MR geometry is exactly the same as a finite linear space with a blocking set that is singled out.
An \emph{isomorphism} between two MR geometries is an isomorphism between the underlying finite linear spaces such that points of the same colour are mapped to points of the same colour.

An MR geometry $(\CS',\chi')$ is a \emph{one-point extension} of the MR geometry $(\CS,\chi)$ if $S$ is a subset of $S'$ with one point less than $S'$, and the restriction of $\chi'$ to $S$ is $\chi$.
Note that in the sequel, when we refer to extensions, we will always mean one-point extensions.
An MR geometry $(\CS,\chi)$ is \emph{minimal} if it is not a one-point extension of any MR geometry, equivalently, if for all $P\in\CS$ such that $\CS\setminus\{P\}$ is non-collinear, the restriction of $\chi$ to $\CS\setminus\{P\}$ is not a proper 2-colouring of $\CS\setminus\{P\}$.
Thus an MR geometry is minimal if we cannot obtain an MR geometry by deleting any single point.

\subsection{The Sylvester-Gallai and Motzkin-Rabin theorems}\label{subsec:sg}
Our motivation for finding small MR geometries arose from questions related to the following theorem (see \cite{PS} or \cite{PS2} for a proof):
\begin{mrtheorem}
If a finite non-collinear set $S$ of points in the real plane \textup{(}affine or projective\textup{)} is $2$-coloured, then there exists a line $\ell$ such that $\ell\cap S$ has size at least two and is monochromatic.
\end{mrtheorem}
In the above terminology this theorem asserts that no MR geometry can be embedded into the real projective plane.
It is an open question whether an MR geometry can be embedded into the complex projective plane.
\begin{problem}
Does there exist an MR geometry that can be embedded into the complex projective plane?
\end{problem}
Since none of the MR geometries on at most $15$ points are embeddable in the projective plane over a division ring of characteristic $\neq 2$, we have the following corollary to Theorem~\ref{mainthm}:
\begin{corollary}
If there exists a finite set $S$ of non-collinear points in the projective plane over a division ring of characteristic $\neq 2$ that can be $2$-coloured such that no line determined by two points in $S$ is monochromatic, then $S$ has size at least $16$.
\end{corollary}

A theorem closely related to the Motzkin-Rabin theorem is the Syl\-ves\-ter-Gallai theorem \cite{MR92b:52010}:
\begin{sgtheorem}
For any finite non-collinear set $S$ of points in the real plane \textup{(}affine or projective\textup{)} there exists a line $\ell$ such that $\ell\cap S$ has size two.
\end{sgtheorem}
In terms of finite linear spaces, this theorem asserts that no proper finite linear space (see \ref{subsec:intro}) can be embedded into the real plane.
These two theorems have proofs that are very closely related (see \cite{PS2}), leading to the question whether the one implies the other, given some reduced set of axioms.
For example, if one assumes the axioms I1, I2, I3 and one of the theorems, does the other follow from it?
First of all it is clear that the Motzkin-Rabin theorem does not imply the Sylvester-Gallai theorem, as the Fano plane $\PG(2,2)$ shows.
The converse implication is the following open question:
\begin{problem}
Does an MR geometry always contain a proper linear space as a subspace?
\end{problem}
Our classification implies that if the answer to this question is negative, then a counterexample would need at least $16$ points.
Indeed, all MR geometries on at most $15$ points that are not themselves proper, contain a Fano plane.
\begin{corollary}
If there exists an MR geometry that does not contain a proper f.l.s., then it must have size at least $16$.
\end{corollary}
For a survey on these and related problems, see \cite{PS2}.

\subsection{Overview of the paper}
In the next section we give a description of all MR geometries on $15$ and fewer points.
In Section~\ref{preparations} we introduce special notation and prove general properties of MR geometries.
The remainder of the paper contains the proof of the classification.
In the proof we distinguish between minimal and non-minimal MR geometries.
The non-minimal geometries on $n$ points may be generated by extending the MR geometries on $n-1$ points.
The results are indicated in the next section with proofs in the Appendix.
It then remains to enumerate the minimal MR geometries.
Since in Section~\ref{preparations} it is shown that there is no MR geometry on $<12$ points, an MR geometry on $12$ points must be minimal.
In Section~\ref{twelvethirteen} we show that there is a unique minimal MR geometry on $12$ points (the punctured projective plane of order $3$), and that there is no minimal MR geometry on $13$ points.
In Section~\ref{sec:fourteen} we show that there are only two minimal MR geometries on $14$ points, and in Section~\ref{fifteen} that there is only one minimal MR geometry on $15$ points, thus completing the proof of the classification (Theorems~\ref{mainthm} and \ref{mainthm2}).

\section{Examples of small MR geometries}\label{examples}
\begin{figure}
\begin{center}
\includegraphics{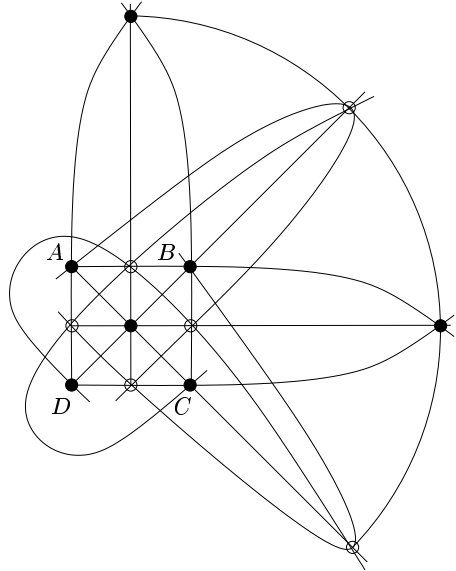}
\end{center}
\caption{$\MR_{13}=\PG(2,3)$---The projective plane of order $3$}\label{fig1}
\end{figure}

We denote the projective plane of order $q=2,3,4$ by $\PG(2,q)$.
See Figure~\ref{fig1} for a standard drawing of $\PG(2,3)$.
It will be convenient for us to represent $\PG(2,4)$ as follows:
The set of points is $\{ i_j : i\in\Z_7, j\in\Z_3\}$, and the lines are
 \[ \{i_j, (i+1)_j, (i+3)_j, f(i)_{j+1}, f(i)_{j+2}\}, \quad (i\in\Z_7, j\in\Z_3),\]
where $f$ is defined by
\[
\begin{tabular}{c|c|c|c|c|c|c|c} 
$i$ & 0 & 1 & 2 & 3 & 4 & 5 & 6 \\ \hline
$f(i)$ & 6 & 3 & 0 & 5 & 1 & 2 & 4\\ 
\end{tabular}.\]
See Figure~\ref{fig2}.
\begin{figure}
\[ \begin{array}{ccc}
0_01_03_06_16_2 & 0_11_13_16_26_0 & 0_21_23_26_06_1 \\
1_02_04_03_13_2 & 1_12_14_13_23_0 & 1_22_24_23_03_1 \\
2_03_05_00_10_2 & 2_13_15_10_20_0 & 2_23_25_20_00_1 \\
3_04_06_05_15_2 & 3_14_16_15_25_0 & 3_24_26_25_05_1 \\
4_05_00_01_11_2 & 4_15_10_11_21_0 & 4_25_20_21_01_1 \\
5_06_01_02_12_2 & 5_16_11_12_22_0 & 5_26_21_22_02_1 \\
6_00_02_04_14_2 & 6_10_12_14_24_0 & 6_20_22_24_04_1
\end{array} \]
\caption{$\PG(2,4)$---The projective plane of order $4$}\label{fig2}
\end{figure}
Clearly each $\{i_j : i\in\Z_7\}$ is a Baer subplane isomorphic to the Fano plane $\PG(2,2)$.

\begin{proposition}\label{exist}
For any $v\geq 12$ there exists an MR geometry on $v$ points.
\end{proposition}

\begin{proof}
Consider a complete quadrilateral in a projective plane $\Proj^2$ of order $q$ with sides $\ell_1,\ell_2,\ell_3,\ell_4$, and set $P_{ij}=\ell_i\cap\ell_j$.
Colour all points on $\ell_1$ and $\ell_2$ green, except for $P_{14}$ and $P_{23}$ which we colour red, and colour all the points on $\ell_3$ and $\ell_4$ red except for $P_{13}$ and $P_{24}$ which we colour green.

We then already have that any line of $\Proj^2$ passes through a green and a red point, except for the line $P_{14}P_{23}$ which does not pass through a green point, and $P_{13}P_{24}$ which does not pass through a red point.
This is rectified by choosing $A\in P_{14}P_{23}$ and $B\in P_{13}P_{24}$ with $A\neq B$, and colouring $A$ green and $B$ red.
This is possible if $q\geq 3$.

We now have a properly $2$-coloured subset of $4q$ points of $\Proj^2$ with the property that every line of $\Proj^2$ already passes through a green and a red point.
We may add points and colour them arbitrarily and still have an MR geometry.
We may thus find an MR geometry for all $v$ satisfying $4q\leq v\leq q^2+q+1$.
By taking $q=2^k$ for $k\geq 3$, we cover all $v\geq 32$.
By taking $q=3,4,5$, we cover all remaining $v\geq 12$ except $14$ and $15$.
These two cases are covered by Examples~\ref{ex:fourteen1} to \ref{ex:fifteen5} below.
\end{proof}

\begin{example}[\boldmath$\MR_{13}=\PG(2,3)$---The projective plane of order $3$\unboldmath]\label{ex:thirteen}\hfill

\par Any complete quadrilateral (see $A,B,C,D$ in Figure~\ref{fig1}) and its three diagonal points ($AB\cap CD, AC\cap BD, AD\cap BC$) form a blocking set.
We thus obtain an MR geometry with six points of one colour and seven of the other.
We denote it by $\MR_{13}$.

It is well-known that $\PG(2,3)$ can be embedded into a Desarguesian projective plane iff the underlying division ring has characteristic $3$.
It is a proper finite linear space.
$\MR_{13}$ is not minimal and it cannot be extended to a $14$-point MR geometry (Lemma~\ref{fourteenminimal}).\qed
\end{example}

\begin{example}[\boldmath$\MR_{12}=\PG(2,3)^\ast$---The punctured projective plane of order $3$\unboldmath]\label{ex:twelve}\hfill
\par Remove any diagonal point of the seven-point blocking set above to obtain a 2-coloured punctured projective plane with six points of each colour, which we denote by $\MR_{12}$.

It is also well-known that the punctured $\PG(2,3)$ can be embedded into a Desarguesian projective plane iff the underlying division ring has characteristic $3$.
It is a proper finite linear space.
By Lemma~\ref{atleastsixpoints}, $\MR_{12}$ is minimal.
It is easily seen that there is only one way to extend $\MR_{12}$ to an MR geometry on $13$ points, viz.\ by putting back the diagonal point.\qed
\end{example}

All the MR geometries on $14$ and $15$ points have subsets isomorphic to the following subset of $\PG(2,4)$:

\begin{lemma}\label{embed}
Let $\D$ be a division ring, and $\CS$ the subspace of $\PG(2,4)$ with points $\{i_0 : i\in\Z_7\}\cup\{2_1,4_1,5_1,6_1\}$.
If $\CS$ can be embedded in $\Proj^2(\D)$, then $\D$ has characteristic $2$ and $x^2+x+1$ has a root in $\D$ \textup{(}equivalently, $\D$ contains the field on four elements\textup{)}, and the embedding of $\CS$ can be extended to an embedding of $\PG(2,4)$ in $\Proj^2(\D)$.

For any two embeddings of\, $\CS$ there is an automorphism of $\Proj^2(\D)$ mapping the one embedded $\PG(2,4)$ onto the other.
\end{lemma}

\begin{proof}
We use projective coordinates.
Without loss of generality we may assume that $0_0=[0,0,1], 1_0=[1,0,1], 5_0=[0,1,1], 3_0=[1,0,0], 4_0=[0,1,0]$ (by applying an appropriate projective transformation of $\Proj^2(\D)$).
See Figure~\ref{fig3}.
\begin{figure}
\begin{center}
\includegraphics{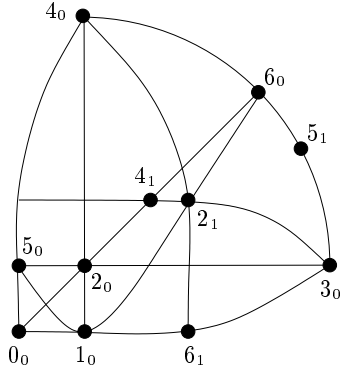}
\end{center}
\caption{Proof of Lemma~\ref{embed}}\label{fig3}
\end{figure}
Then the collinearities $1_02_04_0$ and $2_03_05_0$ force $2_0=[1,1,1]$.
The collinearity $3_04_06_05_1$ implies that $6_0$ and $5_1$ are at infinity.
Also, $0_02_06_0$ forces $6_0=[1,1,0]$ and $1_05_06_0$ forces $6_0=[1,-1,0]$.
It follows that $1=-1$ and the characteristic of $\D$ is $2$.
(Up to now we have merely created the standard coordinatization of the Fano plane $\{i_0 : i\in\Z_7\}$.)

Since $2_1\notin0_03_0$, $2_1\notin0_05_0$, and $2_1\notin3_04_0$, we have $2_1=[a,b,1]$ for some $a,b\neq0$.
By $2_14_06_1$ and $0_03_06_1$ we then have $6_1=[a,0,1]$ (hence $a\neq 1$).
By $0_02_15_1$ and $2_05_16_1$ we have that $0_02_1$ and $2_06_1$ are parallel in the affine plane $\{[x,y,1]:x,y\in\D\}$.
Thus $ba^{-1}=(a+1)^{-1}$, hence $b=(a+1)^{-1}a$.
By $0_02_04_16_0$ we have $4_1=[x,x,1]$ for some $x\neq 0,1$.
By $2_14_13_0$ we then obtain $4_1=[b,b,1]$.
Since $4_1\notin1_04_0$, we have $b\neq 1$.
By $1_04_15_1$ we have that $1_04_1$ is also parallel to $0_02_1$ and $2_06_1$, thus $b(b+1)^{-1}=ba^{-1}$, and $b=a+1$.
Substituting into the previous equation $b=(a+1)^{-1}a$ we obtain $(a+1)^{-1}a=a+1$ from which follows $a^2+a+1=0$.
Now the projective coordinates of all the points of $\CS$ are uniquely determined, as well as the remaining points of $\PG(2,4)$.
\end{proof}

\begin{example}[\boldmath$\MR_{14}^1$\unboldmath---complement of a Baer subplane of \boldmath$\PG(2,4)$\unboldmath]\label{ex:fourteen1}\hfill
\par $\MR_{14}^1$ is the subspace of $\PG(2,4)$ induced by two of the Baer subplanes, say \[\{i_0 : i\in\Z_7\}\cup\{i_1 : i\in\Z_7\}\] to be definite, where we colour each $i_0$ green and each $i_1$ red.

This space is the example of an f.l.s.\ on $14$ points with either $2$ or $4$ points on a line and admitting a blocking set given in \cite{MR97i:05010}.
By Lemma~\ref{embed}, $\MR_{14}^1$ can be embedded into a Desarguesian projective plane iff the underlying division ring has characteristic $2$ and $x^2+x+1$ has a root, and the embedding is projectively unique.
$\MR_{14}^1$ is minimal (proven below in Lemma~\ref{fourteenminimal}).
There are exactly two ways of extending $\MR_{14}^1$ to a $15$-point MR geometry: one is called $\MR_{15}^1$ (Example~\ref{ex:fifteen1}), and the other $\MR_{15}^2$ (Example~\ref{ex:fifteen2}).
See the Appendix for a proof that these are the only one-point extensions.
\qed
\end{example}

\begin{example}[\boldmath$\MR_{14}^2$\unboldmath]\label{ex:fourteen2}\hfill
\par $\MR_{14}^2$ is the subspace of $\PG(2,4)$ induced by \[\{i_0 : i\in\Z_7\}\cup\{2_1,4_1,5_1,6_1,0_2,1_2,3_2\},\] where we colour each $i_0$ green and the remaining points red.
$\MR_{14}^2$ has the same embeddability properties as $\MR_{14}^1$ (Lemma~\ref{embed}), and is minimal (Lemma~\ref{fourteenminimal}).
There are four ways to extend $\MR_{14}^2$: to $\MR_{15}^1, \MR_{15}^{3\mathsf{r}}, \MR_{15}^{3\mathsf{g}}$, and $\MR_{15}^4$ (Examples~\ref{ex:fifteen1}, \ref{ex:fifteen3} and \ref{ex:fifteen4}).
See the Appendix for a proof that there are no other one-point extensions.
\qed
\end{example}

\bigskip
The following four examples on $15$ points are extensions of MR geometries on $14$ points.

\begin{example}[\boldmath$\MR_{15}^1$\unboldmath]\label{ex:fifteen1}\hfill
\par There are two ways to obtain $\MR_{15}^1$ as an extension of an MR geometry on $14$ points.
One is the extension $\MR_{14}^1\cup\{6_2\}$.
The point $6_2$ may be given either colour to give two isomorphic MR geometries.
To be definite we colour $6_2$ green.
The other is the extension $\MR_{14}^2\cup\{6_2\}$, where in this case $6_2$ is forced to be red.

The following is an isomorphism between the two representations:
\begin{center}
\begin{tabular}{c|c|c|c|c|c|c|c|c|c|c|c|c|c|c}
$0_0$ & $1_0$ & $2_0$ & $3_0$ & $4_0$ & $5_0$ & $6_0$ & $0_1$ & $1_1$ & $2_1$ & $3_1$ & $4_1$ & $5_1$ & $6_1$ & $6_2$ \\ \hline
$0_2$ & $3_2$ & $4_1$ & $1_2$ & $2_1$ & $5_1$ & $6_2$ & $1_0$ & $0_0$ & $5_0$ & $3_0$ & $4_0$ & $2_0$ & $6_0$ & $6_1$
\end{tabular}
\end{center}
By Lemma~\ref{embed}, $\MR_{15}^1$ has the same embeddability properties as the MR geometries on $14$ points.\qed
 \end{example}

\begin{example}[\boldmath$\MR_{15}^2$\unboldmath---closed complement of a Baer subplane of \boldmath$\PG(2,4)$\unboldmath]\label{ex:fifteen2}\hfill
\par $\MR_{15}^2$ is the extension of $\MR_{14}^1$ where we add a point $\infty$ to each of the lines $i_0i_1$, $i\in\Z_7$ of $\MR_{14}^1$, and colour it green, say.

It is a proper finite linear space, being case (xii) in Brouwer's classification of proper finite linear spaces on $15$ points \cite{MR83b:05031}.
By Lemma~\ref{embed}, $\MR_{15}^2$ cannot be embedded into a Desarguesian projective plane, since in the (unique) embedding found in the proof of Lemma~\ref{embed}, the lines corresponding to $2_02_1, 4_04_1, 6_06_1$ are not concurrent in $\Proj^2(\D)$.\qed
\end{example}

\begin{example}[\boldmath$\MR_{15}^{3\mathsf{r}}$ and $\MR_{15}^{3\mathsf{g}}$\unboldmath]\label{ex:fifteen3}\hfill
\par Both $\MR_{15}^{3\mathsf{r}}$ and $\MR_{15}^{3\mathsf{g}}$ have the same underlying linear space: $\MR_{14}^2\cup\{4_2\}$.
In $\MR_{15}^{3\mathsf{r}}$ we colour $4_2$ red, and in $\MR_{15}^{3\mathsf{g}}$ we colour it green.

These two MR geometries are not isomorphic (no automorphism of the underlying finite linear space preserves the colour classes).
Indeed, $\MR_{15}^{3\mathsf{r}}$ has a colour class isomorphic to the Fano plane, while no colour class of $\MR_{15}^{3\mathsf{g}}$ is a Fano plane.
By Lemma~\ref{embed}, $\MR_{15}^{3\mathsf{r}}$ and $\MR_{15}^{3\mathsf{g}}$ have the same embeddability properties as the MR geometries on $14$ points.\qed
\end{example}

\begin{example}[\boldmath$\MR_{15}^4$\unboldmath]\label{ex:fifteen4}\hfill
\par $\MR_{15}^4$ is the extension $\MR_{14}^2\cup\{3_1\}$, with $3_1$ coloured red.

By Lemma~\ref{embed}, $\MR_{15}^4$ has the same embeddability properties as the MR geometries on $14$ points.\qed
 \end{example}

\begin{example}[\boldmath$\MR_{15}^5$\unboldmath]\label{ex:fifteen5}\hfill
\par We obtain the underlying f.l.s.\ of $\MR_{15}^5$ by choosing two lines $\ell$ and $m$ of $\PG(2,4)$, and removing four points from $\ell$ and two from $m$, such that the intersection of $\ell$ and $m$ is not among the removed points.
To be definite, we delete $0_1, 1_1, 3_1, 6_2$ from $\ell=6_06_2$ and $0_2, 3_2$ from $m=6_06_1$, and let the green points be $\{0_0, 1_0, 3_0, 4_0, 5_1, 1_2, 5_2\}$, and the red points $\{2_0, 5_0, 6_0, 2_1, 4_1, 6_1, 2_2, 4_2\}$.

This is the unique minimal MR geometry on $15$ points (Proposition~\ref{prop6} in Section~\ref{fifteen}).
$\MR_{15}^5$ is a proper finite linear space, being case (vii) in Brouwer's classification \cite{MR83b:05031}.
By Lemma~\ref{embed}, $\MR_{15}^5$ has the same embeddability properties as the MR geometries on $14$ points.\qed
\end{example}

\begin{figure}
\begin{center}
\includegraphics{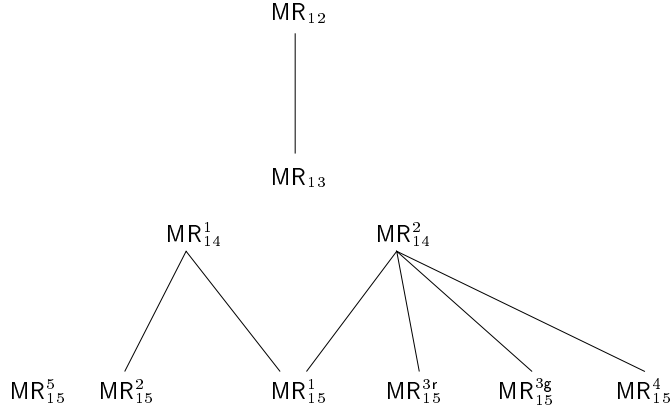}
\end{center}
\caption{The MR geometries on $\leq 15$ points}\label{fig2.5}
\end{figure}

The relationships between the different MR geometries is shown in Figure~\ref{fig2.5}.

We may now state our main result in an exact form.

\begin{theorem}\label{mainthm2}
Up to isomorphism, the only MR geometries on $v\leq 15$ points are
\[\MR_{12}, \MR_{13}, \MR_{14}^1, \MR_{14}^2, \MR_{15}^1, \MR_{15}^2, \MR_{15}^{3\mathsf{r}}, \MR_{15}^{3\mathsf{g}}, \MR_{15}^4, \MR_{15}^5.\]
\end{theorem}

\section{Preparations}\label{preparations}
A \emph{$k$-line} in an f.l.s.\ is a line of size $k$.
The \emph{degree} of a point $P$, denoted by $\deg(P)$, is the number of lines passing through $P$.
A \emph{$d$-point} is a point of degree $d$.

We often draw the \emph{$(P,Q)$-grid} of a finite linear space.
This is the grid formed by the lines through $P$ except $PQ$ drawn vertically, and the lines through $Q$ except $PQ$ drawn horizontally.
The line $PQ$ is then intuitively considered to be at infinity.
In particular, this shows that there are at most $(\deg(P)-1)(\deg(Q)-1)$ points not on $PQ$.

We denote the number of green points in a subset $A$ of an MR geometry by $g(A)$, and the number of red points by $r(A)$.
An \emph{$[\ptype{a}{b}]$-line} in an MR geometry is a line passing through exactly $a$ green and $b$ red points.
We always assume that $a,b\geq 1$.
If we do not want to specify the number of red or green points we put a $\ast$ in its position.
For example, a $[\ptype{\ast}{2}]$-line is a line passing through at least one green point and exactly two red points.
Suppose that the lines passing through $P$ are $[\ptype{g_i}{r_i}]$-lines $(i=1,2,\dots,k)$ where $k=\deg(P)$.
Then the \emph{neighbourhood type} (\emph{ntype} for short) of $P$ is $[\ptype{g_1}{r_1}\ptype{g_2}{r_2}\dots\ptype{g_k}{r_k}]$.
Note that we consider such a matrix to denote the same ntype if the columns are permuted.
We use the same convention regarding $\ast$'s as before.
For example, $[\ptype{\ast}{\ast}\ptype{\ast}{1}\ptype{1}{2}]$ denotes the ntype of a point incident with a $[\ptype{\ast}{\ast}]$-line, a $[\ptype{\ast}{1}]$-line, and a $[\ptype{1}{2}]$-line.

\begin{lemma}\label{neighbourhood}
Let $G$ be a green point in an MR geometry $(\CS,\chi)$ with ntype
\[ \left[\ptype{x}{y}\;\ptype{g_1}{\ast}\;\dots\;\ptype{g_t}{\ast}\;\ptype{1}{\ast}\;\dots\;\ptype{1}{\ast}\right]. \]
Then $ty\leq (t-1)(g(\CS)-x)$.
\textup{(}In particular, $t\geq 2$.\textup{)}
\end{lemma}

\begin{proof}
Denote the $[\ptype{g_i}{\ast}]$-lines by $\ell_i$ $(i=1,\dots,t)$ and the $[\ptype{x}{y}]$-line by $\ell$.
Choose a red point $R_i\in\ell_i$ for each $i=1,\dots,t$.
Let $S_1,\dots,S_y$ be the red points on $\ell$.
Let $g_i'=g_i-1$.
Then $g(\CS)=x+\sum_{i=1}^t g_i'$.
For a fixed $j=1,\dots,t$, the set $R_jS_1\cup R_jS_2\cup\dots \cup R_jS_y$ contains $y$ distinct green points that are also in the set $\bigcup_{\substack{i=1\\i\neq j}}^t\ell_i$.
Therefore, $y\leq\sum_{i\neq j}g_i'$, and
\[ ty \leq \sum_{j=1}^t\sum_{i\neq j}g_i' = (t-1)\sum_{i=1}^t g_i' = (t-1)(g(\CS)-x).\qedhere\]
\end{proof}

The following lemma is very useful in shortening the case analysis when determining the minimal MR geometries of a certain size.
\begin{lemma}\label{minimal}
In a minimal MR geometry $(\CS,\chi)$, each green point $G$ passes through a $[\ptype{1}{b}]$-line for some $b\geq 2$, hence $\deg(G)\leq r(\CS)-1$.
Consequently, each line passes through at most $r(\CS)-1$ points.
\textup{(}An analogous statement holds for red points.\textup{)}
\end{lemma}

\begin{proof}
By minimality, if we remove $G$, then either $\CS\setminus\{G\}$ is collinear, or the 2-colouring is not proper anymore.
In the first case, $\CS$ is a near-pencil (i.e.\ $\CS\setminus\{x\}$ is collinear for some $x$), which is not 2-colourable, a contradiction.
In the second case there is a monochromatic line $\ell\setminus\{G\}$ in $(\CS,\chi)\setminus\{G\}$.
Since $\ell$ is not monochromatic in $(\CS,\chi)$, we must have $G\in\ell$.
Hence all points (at least $2$) of $\ell\setminus\{G\}$ are red.
Since there is a red point on any line through $G$, we then obtain $\deg(G)\leq r(\CS)-1$.

Secondly, given any $[\ptype{a}{b}]$-line $\ell$, there clearly must be a green point $G\notin\ell$.
Then $\deg(G)\geq a+b$.
\end{proof}

The above proof also gives
\begin{lemma}\label{minimal2}
An MR geometry contains a minimal MR geometry.\qed
\end{lemma}

\begin{lemma}\label{atleastfourlines}
The degree of any point in an MR geometry is at least $4$.
\end{lemma}

\begin{proof}
Suppose that a red point $P$ has degree at most $3$.
Let $\ell_1$ and $\ell_2$ be two lines passing through $P$.
Let $G_i$ be a green point on $\ell_i$ $(i=1,2)$.
Let $Q$ be a red point on $G_1G_2$, and let $\ell_3=PQ$.
Then for each red $R \in \ell_2\setminus\{P\}$ there is a green point $G=G_R\in Q R \cap \ell_1$ (Figure~\ref{fig4}).
\begin{figure}
\begin{center}
\includegraphics{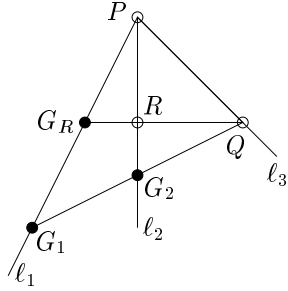}
\end{center}
\caption{Proof of Lemma~\ref{atleastfourlines}}\label{fig4}
\end{figure}
The function $R \mapsto G_R$ is injective.
Thus $r(\ell_2)-1 \leq g(\ell_1\setminus\{G_1\})=g(\ell_1)-1$.
By using a green point on $\ell_3$ we obtain
in a similar way an injection $G \mapsto R_G$ from the green points of
$\ell_1$ to the red points on $\ell_2 \setminus \{P\}$.
Thus $g(\ell_1)\leq r(\ell_2)-1$, a contradiction.
\end{proof}

Although we will not use it, the above lemma can easily be strengthened to show that an MR geometry is not contained in the union of any three lines.

The following lemma shows why an MR geometry has to have at least $12$ points.
\begin{lemma}\label{atleastsixpoints}
In an MR geometry $(\CS,\chi)$,
\begin{enumerate}
\item for any line $\ell$ there are at least $3$ points of each colour not on $\ell$, and
\item there are at least $6$ points of each colour.
\end{enumerate}
\end{lemma}
\begin{proof}
Let $G$ be any green point on $\ell$.
Let $\ell_1,\ell_2,\ell_3$ be three other lines through $G$ (by Lemma~\ref{atleastfourlines}).
Then $\ell_1,\ell_2,\ell_3$ each contains a red point $\notin\ell$.
In a similar way we can find three green points $\notin\ell$.

Next we find six green points.
If we choose two green points and let $\ell$ be the line through them,
by the first part of the proof we obtain $g(\CS)\geq 5$.
If some three green points are collinear, then we would similarly obtain $g(\CS)\geq 6$, and then we are finished.
We may now assume that $g(\CS)=5$ and that no three green points are collinear, and find a contradiction.

By Lemma~\ref{minimal2} we may assume without loss of generality that $(\CS,\chi)$ is minimal.
Thus for every red point $R$ there is an $[\ptype{a}{1}]$-line passing through $R$ with $a\geq 2$, and also $\deg(R)\leq g(\CS)-1=4$ (Lemma~\ref{minimal}).
By Lemma~\ref{atleastfourlines} we then have $\deg(R)=4$ and $a=2$.
Thus $R$ has ntype $[\ptype{1}{\ast}\ptype{1}{\ast}\ptype{1}{\ast}\ptype{2}{1}]$.
It follows that there are exactly as many red points as pairs of green points:
$r(\CS)=\binom{g(\CS)}{2}=10$.

By the pigeon-hole principle there are at least three red points $\neq R$ on one of the $[\ptype{1}{\ast}]$-lines through $R$.
Let $G$ be the green point on this line.
Then $GR$ is a $[\ptype{1}{y}]$-line with $y\geq 4$.
Apply Lemma~\ref{neighbourhood} to $G$ to obtain $4t\leq ty\leq (t-1)(g(\CS)-1)=4(t-1)$, a contradiction.

Therefore, $g(\CS)>5$.
Similarly, $r(\CS)\geq 6$.
\end{proof}

The following generalization of the counting argument for finding the number of red points in the previous proof will be crucial in the determination of the minimal MR geometries.
We first define the weight of a point.
Let $P$ be a point with ntype $[\ptype{g_1}{r_1}\dots\ptype{g_k}{r_k}]$.
Then the \emph{weight} of $P$ is defined as
\[ w(P) = \left\{\begin{array}{ll}\sum_{i=1}^k\binom{r_i}{2}/g_i & \text{if $P$ is green,} \\ \sum_{i=1}^k\binom{g_i}{2}/r_i & \text{if $P$ is red.}\end{array}\right. \]
For each point $P$ in a minimal MR geometry we have $w(P)\geq 1$ by Lemma~\ref{minimal}.

\begin{lemma}\label{weight}
In any MR geometry $(\CS,\chi)$,
\[ \sum_{\textup{$G$ green}} w(G) = \binom{r(\CS)}{2}.\]
\end{lemma}
\begin{proof}
Interchange the order of summation and use axiom I1:
\begin{align*}
\sum_{\text{green $G$}} w(G) & = \sum_{\text{green $G$}}\;\;\sum_{\ell\, \ni\, G} \frac{\tbinom{r(\ell)}{2}}{g(\ell)} \\
& = \sum_\ell\;\,\sum_{\text{green $G\in\ell$}}\frac{\tbinom{r(\ell)}{2}}{g(\ell)} \\
& = \sum_\ell\binom{r(\ell)}{2} = \binom{r(\CS)}{2}.\hfill\qedhere
\end{align*}
\end{proof}

The following is our first application of Lemma~\ref{weight}.

\begin{proposition}\label{sixpoints}
An MR geometry with exactly six green points has at most $13$ points.
\end{proposition}
\begin{proof}
Denote the total number of points by $v$, and the number of red $4$-points by $f$.
Since there are only six green points, for any red point $R$ we have $w(R)\leq 3$ if $\deg(R)=4$, and $w(R)\leq 1$ if $\deg(R)\geq 5$.
By Lemma~\ref{weight} we thus obtain
\begin{equation}\label{eq0}
15=\binom{g}{2}=\sum_{R\text{ red}}w(R)\leq 3f + v-6-f,
\end{equation}
hence
\begin{equation}\label{eq1}
v+2f\geq 21.
\end{equation}

Let $G_1, G_2, G_3$ be non-collinear green points (see Lemma~\ref{atleastsixpoints}), $R_1$ a red point on $G_2G_3$, and $R_2$ a red point on $G_1G_3$.
Then $\deg(R_1), \deg(R_2)\leq 5$.
Let $R_1R_2$ be an $[\ptype{x}{y}]$-line.
Let $G$ be a green point on $R_1R_2$ and $R$ a red point on $GG_1$.
Then $\deg(R)\leq 5$ as well, and considering the $(R_1,R_2)$-grid,
\begin{equation}\label{eq2a}
5\geq \deg(R)\geq x+y\geq v-16,
\end{equation}
and
\begin{equation}\label{eq2}
v\leq 21.
\end{equation}

Suppose $f=0$.
Then by \eqref{eq1} and \eqref{eq2}, $v=21$, by \eqref{eq2a} $R_1R_2$ is a $5$-line, and by \eqref{eq0} the weight of each red point is $1$.
Thus any red point lies on a $[\ptype{2}{1}]$-line, and we may repeat the above argument for any two red points $R_1$ and $R_2$, to obtain that any line passing through at least two red points must be a $5$-line.
It also follows that $\deg(R)=5$ for any red $R$, by choosing $R_1$ and $R_2$ such that $R\notin R_1R_2$.
Thus any red point must have ntype $[\ptype{1}{r_1}\ptype{1}{r_2}\ptype{1}{r_3}\ptype{1}{r_4}\ptype{2}{1}]$ with each $r_i\in\{1,4\}$.
Then $v\leq 19$, a contradiction.

Thus $f\geq 1$.
Let $R_1$ be a red $4$-point, $G_1$ and $G_2$ two red points not collinear with $R_1$ (see Lemma~\ref{atleastsixpoints}), $R_2$ a red point on $G_1G_2$, and $R_1R_2$ an $[\ptype{x}{y}]$-line.
Then $\deg(R_2)\leq 5$ and as before, by considering the $(R_1,R_2)$-grid we obtain
\[ 5\geq x+y\geq v-12.\]
Then $v\leq 17$, and by \eqref{eq1}, $f\geq 2$.
Now let $R_3$ be a second red $4$-point, and $R_1R_3$ an $[\ptype{a}{b}]$-line.
By considering the $(R_1,R_3)$-grid, $5\geq a+b\geq v-9$, and $v\leq 14$.

Suppose $v=14$.
Then we must have $a+b=5$, and by \eqref{eq1}, $f\geq 4$.
Since each red point not on the $5$-line $R_1R_3$ must have degree at least $5$, we have $f\leq b\leq 4$.
Thus $f=4$, $R_1R_3$ is a $[\ptype{1}{4}]$-line, and all red points on $R_1R_3$ are $4$-points.
By Lemma~\ref{neighbourhood} applied to the green point $G$ on $R_1R_3$ we obtain for $t=\deg(G)-1$ that $4t\leq (t-1)(6-1)$, hence $t\geq 5$.
This gives at least $4+5$ red points, a contradiction.

Thus $v\leq 13$.
\end{proof}

It is also possible to prove that an f.l.s.\ with a blocking set of $7$ points has at most $21$ points.
The projective planes of order $3$ (Example~\ref{ex:thirteen}) and $4$ show that these two results are sharp.
In general, using the techniques of Proposition~\ref{sixpoints}, we only get that an f.l.s.\ with a blocking set of $g$ points has at most $g^2-cg$ points.
As the Desarguesian projective planes of square order show, there exist finite linear spaces with a blocking set of size $g=q+\sqrt{q}+1$ and with $q^2+q+1=g^2-g\sqrt{4g-3}+1$ points.

\section{\!Determining the MR geometries on at most \boldmath$13$\unboldmath\ points}\label{twelvethirteen}
\begin{proposition}\label{twelveunique}
There is a unique MR geometry on $12$ points.
\end{proposition}
\begin{proof}
Let $(\CS,\chi)$ be an MR geometry on $12$ points.
By Lemma~\ref{atleastsixpoints}, $g(\CS)=r(\CS)=6$ and $\CS$ is minimal.
By Lemmas~\ref{minimal} and \ref{atleastfourlines} we then have that $4\leq \deg(P)\leq 5$ for all $P\in\CS$.
It follows that each green point $G$ has one of the following ntypes:
\begin{description}
\item[Type 1] $[\ptype{\ast}{1}\ptype{\ast}{1}\ptype{\ast}{1}\ptype{\ast}{1}\ptype{1}{2}]$: $w(G)=1$.
\item[Type 2] $[\ptype{\ast}{1}\ptype{\ast}{1}\ptype{\ast}{1}\ptype{1}{3}]$: $w(G)=3$.
\item[Type 3] $[\ptype{\ast}{1}\ptype{\ast}{1}\ptype{a}{2}\ptype{1}{2}]$ for some $a\geq 1$:
By Lemma~\ref{atleastsixpoints} there are at least $3$ green points not on the $[\ptype{a}{2}]$-line, hence $1\leq a\leq 3$.
There must be a green point on the line through the red points on the two $[\ptype{\ast}{1}]$-lines, giving $a\geq 2$.
By Lemma~\ref{neighbourhood}, $a=3$ is impossible.
Thus $a=2$ and $w(G)=3/2$.
\end{description}
By Lemma~\ref{weight} we have $\sum_{\text{$G$ green}}w(G)=\binom{6}{2}=15$.
By letting $x_i$ be the number of green points of type $i$, this equation becomes
$1x_1+3x_2+\frac32x_3=15$ where we also have $x_1+x_2+x_3=6, x_i\geq 0$.
The only integral solution is $x_1=0, x_2=4, x_3=2$.
Thus four of the green points are of type 2, and two of type 3.
In particular, each green point is a $4$-point.
In the same way, each red point is a $4$-point.
By Vanstone's lemma \cite[Proposition~2.2.3]{MR94m:51019} $\CS$ can be embedded into $\PG(2,3)$, i.e., $\CS$ is the punctured projective plane of order $3$.

It remains to show uniqueness of the colouring up to isomorphism.
Since each line of the punctured projective plane passes through at least three points, and each point is of type 2 or 3, each line must be a $[\ptype{1}{2}]$-, $[\ptype{2}{1}]$-, $[\ptype{2}{2}]$-, $[\ptype{1}{3}]$- or $[\ptype{3}{1}]$-line.
Consider the two green points $G_1$ and $G_2$ of type 3:
Each has ntype $[\ptype{3}{1}\ptype{3}{1}\ptype{2}{2}\ptype{1}{2}]$, $G_1G_2$ is a $[\ptype{2}{2}]$-line, and we may label the remaining points as in Figure~\ref{fig5}(a).
\begin{figure}
\begin{center}
\includegraphics{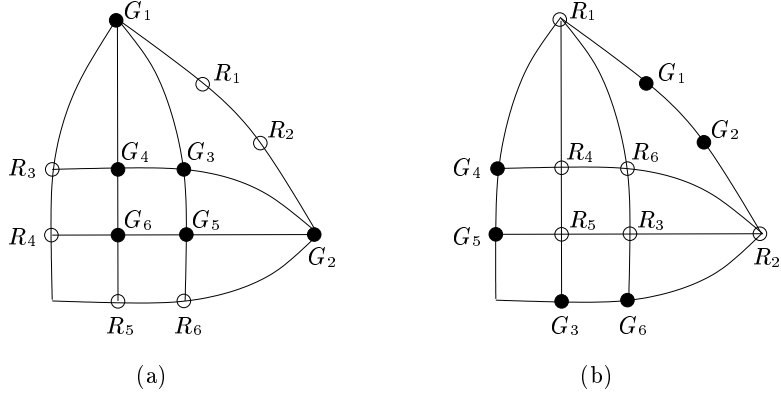}
\end{center}
\caption{Proof of Proposition~\ref{twelveunique}}\label{fig5}
\end{figure}
Since $G_3$ is of type 2, it is on a $[\ptype{1}{3}]$-line.
Thus we have the collinearity $G_3R_4R_5$, and similarly, $G_4R_4R_6$, $G_5R_3R_5$ and $G_6R_3R_6$.
Also we have either $G_3R_4R_1$ or $G_3R_4R_2$.
Without loss of generality we may assume $G_3R_4R_1$ (otherwise relabel $R_1\leftrightarrow R_2$).
Then $G_4R_2R_4$ is also fixed.

In a similar way we consider the red points.
Since there is only one $[\ptype{2}{2}]$-line, we know that $R_1$ and $R_2$ are the two red points playing a similar role to $G_1$ and $G_2$.
We then obtain the two further collinearities $G_5R_2R_3R_5$ and $G_6R_1R_3R_6$ (Figure~\ref{fig5}(b)).
Now we have determined all lines containing at least two green points, and all lines containing at least two red points, hence all lines.
Therefore, we have uniqueness up to isomorphism of MR geometries.
\end{proof}

\begin{proposition}\label{prop:thirteen}
There is a unique MR geometry on $13$ points.
\end{proposition}
\begin{proof}
If $(\CS,\chi)$ is not minimal, it is an extension of the $12$-point MR geometry, which can only be $\PG(2,3)$ (e.g.\ by Vanstone's lemma), and then there are six points of one colour and seven of the other.

We now show that there is no minimal MR geometry on $13$ points.
Suppose there is.
By Lemma~\ref{atleastsixpoints}, $g(\CS)=6$ and $r(\CS)=7$ without loss of generality.

\begin{claim}\label{claim1}
The weight of a green point $G$ is $\leq 7/2$, with equality iff the ntype of $G$ is $[\ptype{3}{1}\ptype{3}{1}\ptype{2}{2}\ptype{1}{3}]$.
\end{claim}
\begin{proof}[Proof of Claim~\ref{claim1}]
By Lemma~\ref{minimal}, $\deg(G)\leq 6$.
If $\deg(G)=6$ then clearly $w(G)=1$, and if $\deg(G)=5$ then clearly $w(G)\leq 3$.
By Lemma~\ref{atleastfourlines}, the only remaining case is $\deg(G)=4$.
Let $\ell$ be a $[\ptype{1}{y}]$-line through $G$ where $y\geq 2$ (Lemma~\ref{minimal}).
By Lemma~\ref{neighbourhood}, $y\leq 3$.

If $y=2$, then $G$ has ntype $[\ptype{\ast}{1}\ptype{\ast}{1}\ptype{a}{3}\ptype{1}{2}]$ or $[\ptype{\ast}{1}\ptype{\ast}{2}\ptype{a}{2}\ptype{1}{2}]$ for some $a\geq 1$.
In the first case, since there must be a green point on the line through the red points on the two $[\ptype{\ast}{1}]$-lines, we must have $a\geq 2$, hence $w(G)\leq5/2$.
In the second case we already have $w(G)\leq 3$.

If $y=3$, then the ntype of $G$ must be $[\ptype{\ast}{1}\ptype{\ast}{1}\ptype{a}{2}\ptype{1}{3}]$, where $a\geq 2$ as before.
Thus $w(G)\leq 7/2$.
The case of equality implies $a=2$.
By considering the lines through the red point on one $[\ptype{\ast}{1}]$-line and the two red points on the $[\ptype{2}{2}]$-line, we obtain that the other $[\ptype{\ast}{1}]$-line has at least three green points.
Hence each $[\ptype{\ast}{1}]$-line has exactly three green points.
\end{proof}

By Lemma~\ref{weight} and Claim~\ref{claim1} we have $21=6\times\frac72\geq\sum_{\text{$G$ green}}w(G)=\binom{r(\CS)}{2}=21$.
Thus we obtain that the ntype of each green point is $[\ptype{3}{1}\ptype{3}{1}\ptype{2}{2}\ptype{1}{3}]$.
Therefore, each line is a $[\ptype{3}{1}]$-, $[\ptype{2}{2}]$-, or $[\ptype{1}{3}]$-line.
By minimality (Lemma~\ref{minimal}), each red point $R$ is on an $[\ptype{a}{1}]$-line where $a\geq 2$.
Thus $a=3$.
Since $\deg(R)\geq 4$ (Lemma~\ref{atleastfourlines}), it follows that $R$ must have ntype $[\ptype{1}{\ast}\ptype{1}{\ast}\ptype{1}{\ast}\ptype{3}{1}]$, hence $[\ptype{1}{3}\ptype{1}{3}\ptype{1}{3}\ptype{3}{1}]$, giving $g(\CS)=6$, a contradiction.

The proof of the uniqueness of the $12$-point MR geometry above can easily be adapted to show that the colouring of $\PG(2,3)$ is unique.
\end{proof}

\section{Determining the MR geometries on \\ \boldmath$14$\unboldmath\ points}\label{sec:fourteen}
\begin{lemma}\label{fourteenminimal}
An MR geometry $(\CS,\chi)$ on $14$ points is minimal.
\end{lemma}
\begin{proof}
Suppose $(\CS,\chi)$ is not minimal.
Then there exists a $P\in\CS$ such that $(\CS,\chi)\setminus\{P\}$ is an MR geometry on $13$ points.
By Proposition~\ref{prop:thirteen} and Example~\ref{ex:thirteen}, $\CS\setminus\{P\}$ is $\PG(2,3)$.
Without loss of generality $P$ is green.

If all lines of $\CS$ through $P$ are $2$-lines, then all the points of $\CS\setminus\{P\}$ must be red.
Thus there is a monochromatic line in $(\CS,\chi)\setminus\{P\}$, a contradiction.
Otherwise there exists a line $\ell$ of $\CS$ through $P$ with at least two more points of $\CS$ on $\ell$.
By Lemma~\ref{atleastsixpoints} there must be a green point $Q\notin\ell$.
Then $PQ$ is a line with at least three points (since it also passes through a red point).
But now $\ell\setminus\{P\}$ and $PQ\setminus\{P\}$ are two non-intersecting lines of $\CS\setminus\{P\}$, contradicting the properties of $\PG(2,3)$.
\end{proof}

\begin{proposition}\label{fourteen}
The only MR geometries on $14$ points are $\MR_{14}^1$ and $\MR_{14}^2$.
\end{proposition}
\begin{proof}
Let $(\CS,\chi)$ be an MR geometry on $14$ points.
By Lemma~\ref{fourteenminimal} it is minimal.
By Lemma~\ref{atleastsixpoints} and Proposition~\ref{sixpoints}, $g(\CS)=r(\CS)=7$.
Without loss of generality we now have two cases:
Either each line of $(\CS,\chi)$ passes through at most three points of either colour, or $(\CS,\chi)$ has a line passing through at least four red points.

\medskip\noindent
\textbf{Case I. \quad \boldmath$r(\ell), g(\ell)\leq 3$ for each line $\ell$ of $\CS$.\unboldmath}

\begin{claim}\label{claim3}
Any line $\ell$ with $r(\ell)=3$ is a $[\ptype{1}{3}]$-line, and any line $\ell$ with $g(\ell)=3$ is a $[\ptype{3}{1}]$-line.
\end{claim}
\begin{proof}[Proof of Claim~\ref{claim3}]
Suppose that $r(\ell)=3$ and $g(\ell)=a\geq 2$.
Let $G$ be a green point on $\ell$.
By Lemmas~\ref{minimal} and \ref{atleastfourlines} it follows that $G$ has ntype $[\ptype{\ast}{1}\ptype{\ast}{1}\ptype{a}{3}\ptype{1}{2}]$.
Let $R_1, R_2, R_3$ be the red points on $\ell$, and $R$ the red point on one of the $[\ptype{\ast}{1}]$-lines.
Then each $RR_i$ has to intersect the other $[\ptype{\ast}{1}]$-line in a green point, giving at least four green point on each $[\ptype{\ast}{1}]$-line, a contradiction.
\end{proof}

\begin{claim}\label{claim4}
The green point $G$ on any $[\ptype{1}{3}]$-line has ntype $[\ptype{\ast}{1}\ptype{\ast}{1}\ptype{\ast}{1}\ptype{\ast}{1}\ptype{1}{3}]$ and weight $3$.
\end{claim}
\begin{proof}[Proof of Claim~\ref{claim4}]
Clearly $\deg(G)\leq 5$.

If $\deg(G)=4$, then the ntype is $[\ptype{\ast}{1}\ptype{\ast}{1}\ptype{\ast}{2}\ptype{1}{3}]$.
Since each line has at most three green points, the two $[\ptype{\ast}{1}]$-lines each contains at most two green points apart from $G$.
Then the $[\ptype{\ast}{2}]$-line passes through a least three green points, hence exactly three by assumption, contradicting Claim~\ref{claim3}.
Thus there are at most six green points, a contradiction.

Thus $\deg(G)=5$ from which the remainder of the claim follows.
\end{proof}

\begin{claim}\label{claim5}
For any green point $G$, $w(G)\leq 3$ with equality iff $G$ is on a $[\ptype{1}{3}]$-line.
\end{claim}
\begin{proof}[Proof of Claim~\ref{claim5}]
Suppose $G$ is not on a $[\ptype{1}{3}]$-line (otherwise use Claim~\ref{claim4}).
Then by Claim~\ref{claim3} there are at most two red points on each line through $G$.
If $w(G)\geq 3$, then there must be at least three $[\ptype{1}{2}]$-lines, and thus only one further line with a single red point.
But then all green points must be on this fourth line, a contradiction.
Thus $w(G)<3$.
\end{proof}
By Claim~\ref{claim5} and Lemma~\ref{weight} we have $7\times 3\leq\sum_{\text{$G$ green}}w(G)=\binom{7}{2}=21$.
Thus $w(G)=3$ for all green $G$, and by Claims~\ref{claim5} and \ref{claim4} they all have the same ntype $[\ptype{\ast}{1}\ptype{\ast}{1}\ptype{\ast}{1}\ptype{\ast}{1}\ptype{1}{3}]$.
It follows that the red points form a Steiner triple system, hence a Fano plane.
Similarly, the green points form a Fano plane.
It follows that the ntype of each green point is $[\ptype{3}{1}\ptype{3}{1}\ptype{3}{1}\ptype{1}{1}\ptype{1}{3}]$, and of each red point is $[\ptype{1}{3}\ptype{1}{3}\ptype{1}{3}\ptype{1}{1}\ptype{3}{1}]$.
Label the green points $A, B, \dots, G$ and the red points $A', B', \dots, G'$ as in Figure~\ref{fig8}.
\begin{figure}
\begin{center}
\includegraphics{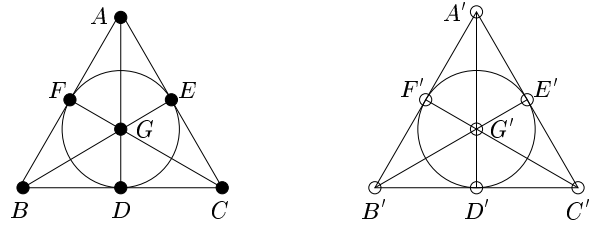}
\end{center}
\caption{The Fano planes in $\MR_{14}^1$ (Proposition~\ref{fourteen})}\label{fig8}
\end{figure}
Since the Fano plane has enough symmetries, without loss of generality $AA'$ is the unique $2$-line passing through $A$, $A'BEG$ the unique $[\ptype{3}{1}]$-line passing through $A'$, and $AB'E'G'$ the unique $[\ptype{1}{3}]$-line passing through $A$.
By permuting $B,E,G$ as well as $B',E',G'$ we obtain without loss of generality the situation in Figure~\ref{fig9}.
\begin{figure}
\begin{center}
\includegraphics{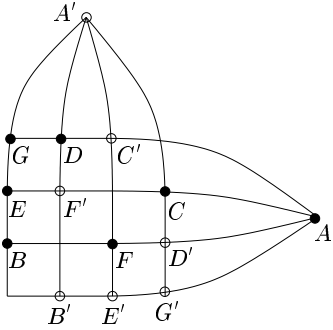}
\end{center}
\caption{Proof of Proposition~\ref{fourteen}}\label{fig9}
\end{figure}
At this stage green and red are symmetric.
We already have the $[\ptype{1}{3}]$-line and one of the $[\ptype{3}{1}]$-lines through $D$.
The remaining two $[\ptype{3}{1}]$-lines through $D$ must then be $DB$ and $DE$.
The only possibility for $DB$ is $BCDE'$ and for $DE$ is $DEFG'$.
It follows that the $[\ptype{1}{1}]$-line through $D$ must be $DD'$.
In the same way the lines through each of $C, F$, and by symmetry, through $C', D', F'$, may be determined.
Finally, to determine the lines through each of $B, E, G, B', E', G'$ is then straightforward.
Thus there is a unique MR geometry in Case I, which is clearly isomorphic to $\MR_{14}^1$.

\medskip\noindent
\textbf{Case II.\quad Line \boldmath$\ell_0$ passes through at least $4$ red points.\unboldmath}

\begin{claim}\label{claim6}
$\ell_0$ is a $[\ptype{1}{4}]$-line, and the ntype of the green point $G$ on $\ell_0$ is $[\ptype{3}{1}\ptype{3}{1}\ptype{3}{1}\ptype{1}{4}]$.
\end{claim}
\begin{proof}[Proof of Claim~\ref{claim6}]
By Lemma~\ref{atleastsixpoints}, $r(\ell_0)=4$.
Fix a green point  $G\in\ell_0$.
Then by Lemma~\ref{atleastfourlines} $\deg(G)=4$, and by Lemma~\ref{minimal} there is no other green point on $\ell_0$.

Let $R_1,R_2,R_3,R_4$ be the red points on $\ell_0$ and $S_1,S_2,S_3$ the red points $\notin\ell_0$.
There is a green point on each of $S_1R_i$ $(i=1,2,3,4)$ hence at least $4$ green points $\neq G$ on $GS_2\cup GS_3$.
Similarly, there are at least $4$ green points $\neq G$ on each $GS_i\cup GS_j$ $(1\leq i<j\leq 3)$.
Since there are $6$ green points $\neq G$, it follows that each $GS_i$ is a $[\ptype{3}{1}]$-line.
\end{proof}

Label the green points $\neq G$ as $G_1, G_2, G_3, H_1, H_2, H_3$ such that we have the lines $GG_iH_iS_i$ $(i=1,2,3)$.

\begin{claim}\label{claim7}
The ntype of each $S_i$ is $[\ptype{1}{2}\ptype{1}{2}\ptype{1}{3}\ptype{1}{3}\ptype{3}{1}]$ \textup{(}in particular, $w(S_i)=3$ and $\deg(S_i)=5$\textup{)}.
Also, each $S_iS_j$ passes through some $R_k$.
\end{claim}
\begin{proof}[Proof of Claim~\ref{claim7}]
Since $S_i$ is not on the $5$-line $\ell_0$, $\deg(S_i)\geq 5$.
Since there are four green points not on $GS_i$, $\deg(S_i)\leq 5$.

Since there must be a green point on each $S_iR_j$ $(j=1,2,3,4)$, we must have the situation in Figure~\ref{fig10}.
\begin{figure}[t]
\begin{center}
\includegraphics{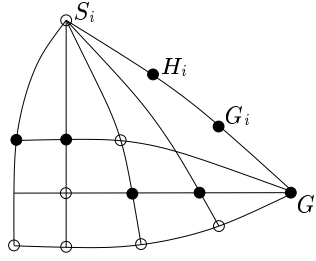}
\end{center}
\caption{Neighbourhood of $S_i$ (Claim~\ref{claim7})}\label{fig10}
\end{figure}
\end{proof}

We now apply Claim~\ref{claim7} to each $S_i$.
For $S_1$ we have without loss of generality the situation in Figure~\ref{fig11} (after possibly permuting $R_1, R_2, R_3, R_4$ and interchanging $S_2\leftrightarrow S_3$, $G_2\leftrightarrow H_2$, $G_3\leftrightarrow H_3$).
\begin{figure}[b]
\begin{center}
\includegraphics{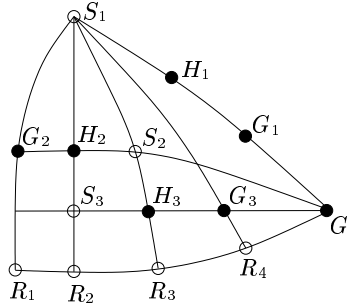}
\end{center}
\caption{Neighbourhood of $S_1$}\label{fig11}
\end{figure}
Thus we have determined the lines $S_1R_1G_2$, $S_1R_4G_3$, $S_1S_2R_3H_3$, $S_1S_3R_2H_2$.
By Claim~\ref{claim7}, $S_2S_3$ passes through some $R_i$.
We must have $i\neq2,3$.
Without loss of generality $S_2S_3$ passes through $R_1$ (if $R_4$, interchange $G_2\leftrightarrow G_3$, $H_2\leftrightarrow H_3$, $R_1\leftrightarrow R_4$, $R_2\leftrightarrow R_3$, $S_2\leftrightarrow S_3$).

Applying Claim~\ref{claim7} to $S_2$ and using the collinearities obtained so far we obtain the situation in Figure~\ref{fig12}.
\begin{figure}
\begin{minipage}[t]{5.6cm}
\begin{center}
\includegraphics{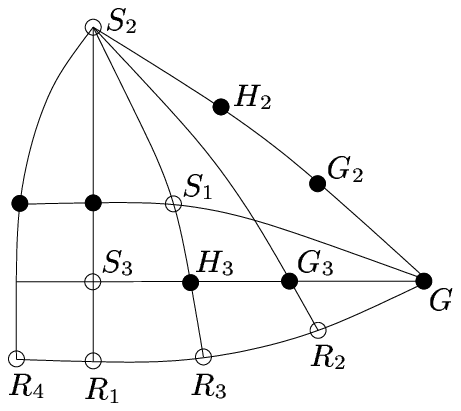}
\end{center}
\caption{Neighbourhood of $S_2$}\label{fig12}
\end{minipage}
\hfill
\begin{minipage}[t]{5.6cm}
\begin{center}
\includegraphics{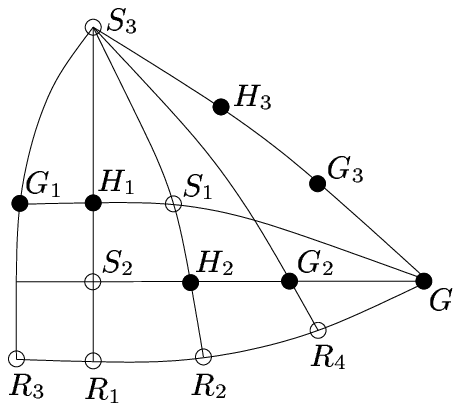}
\end{center}
\caption{Neighbourhood of $S_3$}\label{fig13}
\end{minipage}
\end{figure}
The two green points $GS_1\cap S_2R_4$ and $GS_1\cap S_2R_1$ are $G_1$ and $H_1$.
Without loss of generality, $GS_1\cap S_2R_4=G_1$, $GS_1\cap S_2R_1=H_1$, since we may still interchange $G_1\leftrightarrow H_1$.
Thus we have the additional lines $S_2G_1R_4$, $S_2S_3R_1H_1$, $S_2G_3R_2$.

Applying Claim~\ref{claim7} and the collinearities obtained so far to $S_3$ we obtain the situation in Figure~\ref{fig13}.
In particular we have the lines $S_3R_3G_1$ and $S_3R_4G_2$.

Thus $R_1, R_2, R_3$ each lies on a $[\ptype{1}{4}]$-, $[\ptype{1}{3}]$- and $[\ptype{1}{2}]$-line, and $R_4$ lies on a $[\ptype{1}{4}]$-line and three $[\ptype{1}{2}]$-lines.
By Lemma~\ref{minimal}, each $R_i$ is furthermore on an $[\ptype{a}{1}]$-line, where $a\geq 2$.

\begin{claim}\label{claim8}
For each $i=1,2,3,4$, $w(R_i)\leq3$ with equality iff the ntype of $R_i$ is $[\ptype{1}{\ast}\ptype{1}{\ast}\ptype{1}{\ast}\ptype{1}{\ast}\ptype{3}{1}]$.
\end{claim}
\begin{proof}[Proof of Claim~\ref{claim8}]
If $\deg(R_i)=4$, then $a=4$, hence there is a $[\ptype{4}{1}]$-line not passing through $G$.
Thus $\deg(G)\geq 5$, a contradiction.
Thus $\deg(R_i)\geq 5$, and the claim is now clear.
\end{proof}

We now have $7\times3\geq\sum_{\text{$R$ red}}w(R)=\binom{g(\CS)}{2}=21$, and it follows that each $R_i$ has ntype $[\ptype{1}{\ast}\ptype{1}{\ast}\ptype{1}{\ast}\ptype{1}{\ast}\ptype{3}{1}]$.
By also considering the ntypes of the $S_i$'s (see Claim~\ref{claim7}), we see that each line passes through exactly $1$ or $3$ green points.
Thus the green points form a Steiner triple system, hence a Fano plane.
We already have $R_1R_2R_3R_4G$, $R_4S_1G_3$, $R_4S_2G_1$, $R_4S_3G_2$ as lines through $R_4$.
Thus the fifth line must be $R_4H_1H_2H_3$.
We now have $GG_iH_i$ $(i=1,2,3)$ and $H_1H_2H_3$ as collinear green triples.
This determines the Fano plane:
the remaining green collinear triples must be $H_iG_{i+1}G_{i+2}$ $(i=1,2,3)$ (subscripts modulo $3$).

Considering the neighbourhoods of $R_1, R_2, R_3$ and the green collinearities determined above we obtain the remaining lines of $(\CS,\chi)$: $H_iR_{i+1}$, $R_iG_iG_{i+2}H_{i+1}$ $(i=1,2,3)$ (subscripts mod $3$).
Thus we have found a unique linear space with a unique colouring in Case II, which is clearly $\MR_{14}^2$.
\end{proof}

\section{Determining the MR geometries on\\ \boldmath$15$\unboldmath\ points}\label{fifteen}
The proof of the following lemma is routine, for example, by checking the list of finite linear spaces on $7$ points in \cite[Appendix]{MR94m:51019}.

\begin{lemma}\label{sevenpointfls}
There is only one finite linear space on $7$ points with six $3$-lines and three $2$-lines, namely the Fano plane with one line $ABC$ replaced by the three lines $AB$, $BC$, $AC$.
See Figure~\ref{fig14}.
\begin{figure}
\begin{center}
\includegraphics[scale=0.45]{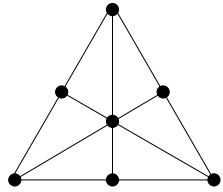}
\end{center}
\caption{The unique f.l.s.\ on $7$ points with six $3$-lines and three $2$-lines}\label{fig14}
\end{figure}
\end{lemma}

The non-minimal MR geometries on $15$ points are enumerated in Examples~\ref{ex:fifteen1} to \ref{ex:fifteen4} and the Appendix.
\begin{proposition}\label{prop6}
There is exactly one minimal MR geometry on $15$ points.
\end{proposition}
\begin{proof}
By Lemmas~\ref{atleastsixpoints} and \ref{sixpoints} there are at least $7$ green points and at least $7$ red points.
Without loss of generality assume $g(\CS)=7$ and $r(\CS)=8$.
\begin{claim}\label{claimzero}
There is no $[\ptype{a}{b}]$-line if $a\geq 3, b\geq 2$.
\end{claim}
\begin{proof}[Proof of Claim~\ref{claimzero}]
Let two of the red points on an $[\ptype{a}{b}]$-line be $R_1$ and $R_2$.
By Lemma~\ref{minimal} there is also a $[\ptype{c}{1}]$-line through $R_1$ for some $c\geq 2$.
Since there are $a+c\geq 5$ green points on these two lines, there are at most two further green points, hence at most two further lines through $R_1$.
Thus $\deg(R_1)\leq 4$, and by Lemma~\ref{atleastfourlines}, $\deg(R_1)=4$, and $a=3, c=2$.
It follows that the ntype of $R_1$, and similarly, of $R_2$, must be $[\ptype{1}{\ast}\ptype{1}{\ast}\ptype{3}{b}\ptype{2}{1}]$.
By considering the $(R_1,R_2)$-grid, we see that there are in total at most $3^2-1 + 3+b$ points, hence $b\geq 4$.
Thus there are at least seven points on $R_1R_2$, contradicting Lemma~\ref{minimal}.
\end{proof}

\begin{claim}\label{claimone}
There is no green $4$-point.
\end{claim}
\begin{proof}[Proof of Claim~\ref{claimone}]
Suppose $G$ is a green $4$-point.
Let the lines through $G$ be $GR_i$, $i=1,2,3,4$, and let $GR_i$ be an $[\ptype{x_i}{y_i}]$-line.
By Lemma~\ref{minimal} we may assume without loss that $x_4=1, y_4\geq 2$.

Suppose that $x_1\geq 4$.
Then the number of green points not on $GR_1$ is at most three.
Since for each red $R\in GR_4$ there is a different green point on $R_1R$, we then must have $y_4\leq 3$.
Also, by Claim~\ref{claimzero}, $y_1=1$.
Then $y_2+y_3\geq 4$.
By Lemma~\ref{atleastsixpoints} there are at least three green points not on $GR_1$, hence $x_2\geq 3$ or $x_3\geq 3$.
Without loss we consider only the case $x_2\geq 3$.
By Claim~\ref{claimzero}, $y_2=1$.
Thus $y_3\geq 3$.
For each red $R\in GR_3$ there is a distinct green point on $R_1R$.
These green points must all lie on $GR_2$, thus $x_2\geq 4$.
Since $x_1\geq 4$ as well, all green points now lie on $GR_1\cup GR_2$.
Thus there is no green point on $R_1R_2$, a contradiction.

Therefore, $x_1\leq 3$, and similarly, $x_2,x_3\leq 3$.
Since there are only seven green points, we must have $x_1=x_2=x_3=3$.
By Claim~\ref{claimzero}, $y_1=y_2=y_3=1$.
Thus $y_4=5$.
For each red $R\in GR_4$ there is a different green point on $R_1R$.
These five green points must all lie on $GR_2\cup GR_3$.
Therefore, $x_2+x_3\geq 7$, a contradiction.
\end{proof}

\begin{claim}\label{claimtwo}
There is at most one red $4$-point.
\end{claim}
\begin{proof}[Proof of Claim~\ref{claimtwo}]
Suppose $R_1$ and $R_2$ are two red $4$-points.
Let $R_1R_2$ be an $[\ptype{x}{y}]$-line.
By considering the $(R_1,R_2)$-grid, we see that there are at least $15-9=6$ points on $R_1R_2$.
Thus $x+y\geq 6$.
By Lemma~\ref{minimal}, $x+y=6$ and the $3\times3$ grid is full.
Again by Lemma~\ref{minimal} both $R_1$ and $R_2$ lie on $[\ptype{3}{1}]$-lines.
Thus there are at least five green points not on $R_1R_2$, and $x\leq 2$.

If $x=1$, then $y=5$, and there are three red points not on $R_1R_2$.
Thus the degree of the green point on $R_1R_2$ is at most $4$, contradicting Claim~\ref{claimone}.

If $x=2$, then $y=4$, and by Lemma~\ref{minimal}, through any green $G\in R_1R_2$ there must be some $[\ptype{1}{a}]$-line, $a\geq 2$.
Then $\deg(G)\leq 4$, again a contradiction.
\end{proof}

\begin{claim}\label{claimthree}
For any green point $G$ exactly one of the following cases holds:
\begin{enumerate}
\item $w(G)=6$ and the ntype of $G$ is $[\ptype{\ast}{1}\ptype{\ast}{1}\ptype{\ast}{1}\ptype{\ast}{1}\ptype{1}{4}]$.
\item $w(G)=4$ and the ntype of $G$ is $[\ptype{\ast}{1}\ptype{\ast}{1}\ptype{\ast}{1}\ptype{1}{2}\ptype{1}{3}]$.
\item $w(G)<4$.
\end{enumerate}
\end{claim}
\begin{proof}[Proof of Claim~\ref{claimthree}]
By Lemma~\ref{minimal}, $\deg(G)\leq 7$, and by Claim~\ref{claimone}, $\deg(G)\geq 5$.

If $\deg(G)=7$, then by Lemma~\ref{minimal} its ntype must be $[\ptype{\ast}{1}\ptype{\ast}{1}\ptype{\ast}{1}\ptype{\ast}{1}\ptype{\ast}{1}\ptype{\ast}{1}\ptype{1}{2}]$, and $w(G)=1$.

If $\deg(G)=6$, then its ntype must be either $[\ptype{\ast}{1}\ptype{\ast}{1}\ptype{\ast}{1}\ptype{\ast}{1}\ptype{\ast}{1}\ptype{1}{3}]$ with $w(G)=3$, or $[\ptype{\ast}{1}\ptype{\ast}{1}\ptype{\ast}{1}\ptype{\ast}{1}\ptype{\ast}{2}\ptype{1}{2}]$ with $w(G)\leq 2$.

If $\deg(G)=5$, its ntype must be one of:
\begin{itemize}
\item $[\ptype{\ast}{1}\ptype{\ast}{1}\ptype{\ast}{1}\ptype{\ast}{1}\ptype{1}{4}]$ with $w(G)=6$,
\item $[\ptype{\ast}{1}\ptype{\ast}{1}\ptype{\ast}{1}\ptype{\ast}{2}\ptype{\ast}{3}]$ with $w(G)\leq 4$, and equality implies that the ntype must be $[\ptype{\ast}{1}\ptype{\ast}{1}\ptype{\ast}{1}\ptype{1}{2}\ptype{1}{3}]$,
\item $[\ptype{\ast}{1}\ptype{\ast}{1}\ptype{\ast}{2}\ptype{\ast}{2}\ptype{1}{2}]$ with $w(G)\leq 3$.\qedhere
\end{itemize}
\end{proof}

\begin{claim}\label{claimfour}
Some green point has weight $6$.
\end{claim}
\begin{proof}[Proof of Claim~\ref{claimfour}]
Suppose no green point has weight $6$.
Then by Claim~\ref{claimthree} the sum of the green weights is
\[ 28=\binom{r(\CS)}{2}=\sum_{G\text{ green}}w(G)\leq 4\times 7.\]
Thus each green point has weight $4$ and ntype $[\ptype{\ast}{1}\ptype{\ast}{1}\ptype{\ast}{1}\ptype{1}{2}\ptype{1}{3}]$ by Claim~\ref{claimthree}.
In particular, each line has at most $3$ red points, and a line passing through at least $2$ red points passes through only one green point.

We now consider the weights of the red points.
Let all the $[\ptype{1}{\ast}]$-lines through a red point $R$ be $RG_1,\dots,RG_n$.
Then by the above, $RG_i$ is a $[\ptype{1}{b_i}]$-line for some $b_i\leq 3$.
Also, all other lines through $R$ must be $[\ptype{\ast}{1}]$-lines.
Adding up all the red points on the lines through $R$ we find $7=\sum_{i=1}^n(b_i-1)\leq 2n$, hence $n\geq 4$.
By Lemma~\ref{minimal}, there is an $[\ptype{a}{1}]$-line through $R$, with $a\geq 2$.
Let $H_1$ and $H_2$ be two green points on this line.
We now already have six of the green points: $G_1,\dots,G_4,H_1,H_2$.
The seventh green point can be either on $H_1H_2$, giving $w(R)=3$, or on a sixth line through $R$, giving $w(R)=1$.

Thus the sum of the red weights is a sum of eight odd numbers ($1$ or $3$), which is even.
By Lemma~\ref{weight}, this sum must equal $\binom{g(\CS)}{2}=21$, a contradiction.
\end{proof}

Fix a green point $G$ of weight $6$.
By Claim~\ref{claimthree}, $G$ has ntype $[\ptype{\ast}{1}\ptype{\ast}{1}\ptype{\ast}{1}\ptype{\ast}{1}\ptype{1}{4}]$.
Denote the red points on the $[\ptype{1}{4}]$-line through $G$ by $R_1$, $R_2$, $R_3$, $R_4$.
Thus we have the full line \[ GR_1R_2R_3R_4.\tag{1}\]

\begin{claim}\label{claimfive}
$G$ has ntype either
\begin{description}
\item[{\upshape \bfseries Case \boldmath$\alpha$\unboldmath:}] $[\ptype{1}{1}\ptype{3}{1}\ptype{3}{1}\ptype{3}{1}\ptype{1}{4}]$ or
\item[{\upshape \bfseries Case \boldmath$\beta$\unboldmath:}]  $[\ptype{2}{1}\ptype{2}{1}\ptype{3}{1}\ptype{3}{1}\ptype{1}{4}]$.
\end{description}
\end{claim}
\begin{proof}[Proof of Claim~\ref{claimfive}]
Without loss of generality $G$ has ntype $[\ptype{x_1}{1}\ptype{x_2}{1}\ptype{x_3}{1}\ptype{x_4}{1}\ptype{1}{4}]$ with
\[ 1\leq x_1\leq x_2\leq x_3\leq x_4.\]
Let $S_i$ be the red point on the $[\ptype{x_i}{1}]$-line.
Since there is a different green point on each line $S_4R_i$, $i=1,2,3,4$, and these $4$ green points must be in $GS_1\cup GS_2\cup GS_3\setminus\{G\}$, we obtain
\[ x_1+x_2+x_3\geq 7.\]
We also have \[x_1+x_2+x_3+x_4=10.\]
There are only two solutions to this equation under the above two constraints:
\[ (x_1,x_2,x_3,x_4)=(2,2,3,3), (1,3,3,3).\qedhere\]
\end{proof}

We now show that Case $\alpha$ in Claim~\ref{claimfive} leads to a contradiction, while Case $\beta$ leads to exactly one MR geometry.

\subsection*{\boldmath Case $\alpha$: $G$ has ntype $[\ptype{1}{1}\ptype{3}{1}\ptype{3}{1}\ptype{3}{1}\ptype{1}{4}]$\unboldmath}
Let $S$ be the red point on the $[\ptype{1}{1}]$-line through $G$, and let $S_1, S_2, S_3$ be the red points on the three $[\ptype{3}{1}]$-lines through $G$.
We now consider the weights of the red points.
\begin{claim}\label{a1}
\mbox{}\begin{itemize}
\item $S$ has weight $1$ and ntype $[\ptype{1}{\ast}\ptype{1}{\ast}\ptype{1}{\ast}\ptype{1}{\ast}\ptype{1}{\ast}\ptype{2}{1}]$.
\item Each $S_i$ has weight $3$ and ntype $[\ptype{1}{\ast}\ptype{1}{\ast}\ptype{1}{\ast}\ptype{1}{\ast}\ptype{3}{1}]$.
\end{itemize}
\end{claim}
\begin{proof}[Proof of Claim~\ref{a1}]
The lines $SG$ and $SR_i$, $i=1,2,3,4$, give $\deg(S)\geq 5$.
By Lemma~\ref{minimal}, $S$ must also be on a sixth $[\ptype{a}{1}]$-line, $a\geq 2$.
It follows that $a=2$, $\deg(S)=6$, and then $S$ must have the stated ntype and weight.

Without loss consider $S_1$.
The line $S_1G$ passes through three green points, and each $S_1R_i$, $i=1,2,3,4$, passes through at least one green point.
Since there are seven green points, it follows that each $S_1R_i$ passes through exactly one green point, and $\deg(S_1)=5$.
Thus $S_1$ has the stated ntype and weight.
\end{proof}

\begin{claim}\label{a2}
For each $R_i$, $i=1,2,3,4$, we have one of the following:
\begin{itemize}
\item $\deg(R_i)=6$ and $w(R_i)=1$,
\item $\deg(R_i)=5$ and either $w(R_i)<3/2$ or $w(R_i)\in\{3/2,2,3\}$,
\item $\deg(R_i)=4$ and $3<w(R_i)\leq 7/2$.
\end{itemize}
For the following weights, the ntypes must be as in the table:
\begin{center}\renewcommand{\arraystretch}{1.5}
\begin{tabular}{c|c}
$w(R_i)$ & ntype\\ \hline
$3/2$ & $[\ptype{1}{\ast}\ptype{1}{\ast}\ptype{1}{\ast}\ptype{2}{2}\ptype{2}{1}]$ \\
$2$ & $[\ptype{1}{\ast}\ptype{1}{\ast}\ptype{1}{\ast}\ptype{2}{1}\ptype{2}{1}]$ \\
$3$ & $[\ptype{1}{\ast}\ptype{1}{\ast}\ptype{1}{\ast}\ptype{1}{\ast}\ptype{3}{1}]$ \\
$7/2$ & $[\ptype{1}{\ast}\ptype{1}{\ast}\ptype{2}{2}\ptype{3}{1}]$
\end{tabular}
\end{center}
\end{claim}
\begin{proof}[Proof of Claim~\ref{a2}]
By Lemmas~\ref{minimal} and \ref{atleastfourlines} we have $4\leq\deg(R_i)\leq 6$.
The cases $\deg(R_i)=6$ and $\deg(R_i)=5$ are clear.
In the case $\deg(R_i)=4$ we consider the $(R_i,G)$-grid.
Note that $GR_i$ is a $[\ptype{1}{4}]$-line and $GS$ is a $[\ptype{1}{1}]$-line.
Thus the complement of $GR_i\cup GS$ must consist of six green points and three red points.
See Figure~\ref{fig15}.
\begin{figure}
\begin{center}
\begin{overpic}[scale=0.5]{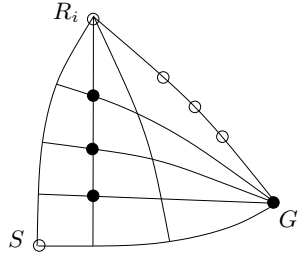}\small
\put(8,94){$R_i$}
\put(99,10){$G$}
\put(-10,2){$S$}
\end{overpic}
\end{center}
\caption{$\deg(R_i)=4$ in Claim~\ref{a2}}\label{fig15}
\end{figure}
By this and Lemma~\ref{minimal}, $R_i$ must lie on some $[\ptype{3}{1}]$-line.
Thus $w(R_i)\geq 3$.
By placing the remaining three green points and three red points we see that the remainder of the claim follows.
\end{proof}

Since by Lemma~\ref{weight} the sum of the red weights must be $\binom{7}{2}$, by Claim~\ref{a1} we have $\sum_{i=1}^4 w(R_i)=11$.
Suppose that at most two $R_i$'s satisfy $w(R_i)\geq 3$.
Then by Claim~\ref{a2}, the remaining $R_i$'s satisfy $w(R_i)\leq 2$, and recalling that there is at most one red $4$-point (Claim~\ref{claimtwo}), we obtain
\[ 11=\sum_{i=1}^4 w(R_i)\leq 7/2+3+2+2,\]
a contradiction.

Thus at least three $R_i$'s have $w(R_i)\geq 3$.
We then have by Claim~\ref{claimtwo} one of the following two cases without loss of generality:
\begin{description}
\item[Case $\alpha1$:] $w(R_1)=w(R_2)=w(R_3)=3$.
\item[Case $\alpha2$:] $w(R_1)>3, w(R_2)=w(R_3)=3$.
\end{description}

\subsection*{Case \boldmath$\alpha1$\unboldmath}
We must have $w(R_4)=2$.
We now consider the subgeometry of the seven green points.
By examining the ntypes of the red points (see Claims~\ref{a1} and \ref{a2}), we see that there are in total six $[\ptype{3}{\ast}]$-lines and three $[\ptype{2}{\ast}]$-lines.
By Lemma~\ref{sevenpointfls} we see that any two $[\ptype{2}{\ast}]$-lines must intersect in a green point.
However, from the ntype of $R_4$ we see that there are two $[\ptype{2}{\ast}]$-lines not intersecting in a green point, a contradiction.

\subsection*{Case \boldmath$\alpha2$\unboldmath}
We must have $w(R_4)<2$.
By Claim~\ref{a2}, $w(R_4)\leq 3/2$.
Thus $w(R_1)\geq 7/2$.
Again by Claim~\ref{a2}, $w(R_1)=7/2$ and $R_1$ has ntype $[\ptype{1}{4}\ptype{1}{4}\ptype{2}{2}\ptype{3}{1}]$, $w(R_4)=3/2$ and $R_4$ has ntype $[\ptype{1}{\ast}\ptype{1}{\ast}\ptype{1}{\ast}\ptype{2}{2}\ptype{2}{1}]$.
By considering the subgeometry of the seven green points we obtain the same contradiction as in Case $\alpha1$.

\subsection*{\boldmath Case $\beta$: $G$ has ntype $[\ptype{2}{1}\ptype{2}{1}\ptype{3}{1}\ptype{3}{1}\ptype{1}{4}]$\unboldmath}
Let $S_1$ and $S_2$ be the red points on the $[\ptype{2}{1}]$-lines through $G$, and let $T_1$ and $T_2$ be the red points on the $[\ptype{3}{1}]$-lines through $G$.
We again consider the weights of the red points.

\begin{claim}\label{b1}
\mbox{}\begin{itemize}
\item $w(S_i)\in\{1,4/3,3/2\}$.
\item $w(T_i)=3$.
\item If $\deg(R_i)\geq 5$ then $w(R_i)=3$ or $w(R_i)\leq 2$.
\item If $\deg(R_i)=4$ then $w(R_i)\in\{10/3, 7/2, 6\}$.
\end{itemize}
For the following weights the ntypes are as in the table:
\begin{center}\renewcommand{\arraystretch}{1.5}
\hfill\hfill
\begin{tabular}[t]{c|c}
$w(S_i)$ & ntype\\ \hline
$1$ & $[\ptype{1}{\ast}\ptype{1}{\ast}\ptype{1}{\ast}\ptype{1}{\ast}\ptype{1}{\ast}\ptype{2}{1}]$ \\
$4/3$ & $[\ptype{1}{\ast}\ptype{1}{\ast}\ptype{1}{\ast}\ptype{2}{3}\ptype{2}{1}]$ \\
$3/2$ & $[\ptype{1}{\ast}\ptype{1}{\ast}\ptype{1}{\ast}\ptype{2}{2}\ptype{2}{1}]$
\end{tabular}
\hfill
\begin{tabular}[t]{c|c}
$w(T_i)$ & ntype\\ \hline
$3$ & $[\ptype{1}{\ast}\ptype{1}{\ast}\ptype{1}{\ast}\ptype{1}{\ast}\ptype{3}{1}]$
\end{tabular}
\hfill
\begin{tabular}[t]{c|c}
$w(R_i)$ & ntype\\ \hline
$3$ & $[\ptype{1}{\ast}\ptype{1}{\ast}\ptype{1}{\ast}\ptype{1}{\ast}\ptype{3}{1}]$ \\
$10/3$ & $[\ptype{1}{\ast}\ptype{1}{\ast}\ptype{2}{3}\ptype{3}{1}]$ \\
$7/2$ & $[\ptype{1}{\ast}\ptype{1}{\ast}\ptype{2}{2}\ptype{3}{1}]$ \\
$6$ & $[\ptype{1}{\ast}\ptype{1}{\ast}\ptype{1}{\ast}\ptype{4}{1}]$
\end{tabular}
\hfill\hfill\mbox{}
\end{center}
\end{claim}
\begin{proof}[Proof of Claim~\ref{b1}]
Since $GR_1$ is a $5$-line not passing through $S_i$, $\deg(S_i)\geq 5$.
Thus we have either $\deg(S_i)=6$ and $w(S_i)=1$, or $S_i$ has degree $5$ and ntype $[\ptype{1}{\ast}\ptype{1}{\ast}\ptype{1}{\ast}\ptype{2}{a}\ptype{2}{1}]$ for some $a\geq 2$ (consider the lines $S_iG$ and $S_iR_j$), and weight $1+1/a$.
Suppose $a\geq 4$.
Let $H$ be a green point on the $[\ptype{2}{a}]$-line through $S_i$.
By Lemma~\ref{minimal} it follows that $\deg(H)=4$, contradicting Claim~\ref{claimone}.
Thus $a\leq 3$, giving $w(S_i)\in\{4/3,3/2\}$.

It is easily seen that $w(T_i)=3$ and if $\deg(R_i)\geq 5$, then $w(R_i)=3$ or $w(R_i)\leq 2$.

Suppose now $\deg(R_i)=4$.
Consider the $(R_i,T_1)$-grid (Figure~\ref{fig16}).
\begin{figure}
\begin{center}
\begin{overpic}[scale=0.5]{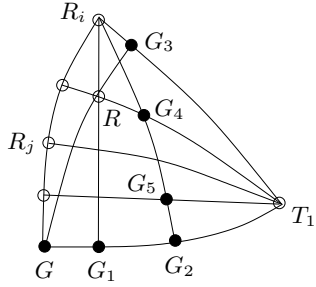}\small
\put(9,94){$R_i$}
\put(102,12){$T_1$}
\put(-1,-10){$G$}
\put(-12,42){$R_j$}
\put(20,-10){$G_1$}
\put(26,51){$R$}
\put(51,-8){$G_2$}
\put(36,25){$G_5$}
\put(47,55){$G_4$}
\put(43,82){$G_3$}
\end{overpic}
\end{center}
\caption{$\deg(R_i)=4$ in Claim~\ref{b1}}\label{fig16}
\end{figure}
Let the green points on $GT_1$ be $G, G_1, G_2$.
Let $G_3$ be the green point on $R_iT_1$.
Let $R$ be a red point on $GG_3$.
Let $G_4$ be the green point on $RT_1$.
By Claim~\ref{claimzero}, there is at most one green point $\neq G_1$ on $G_1R_i$.
Thus there must be a third green point $G_5$ on $R_iG_2$.
Let $R_j$ be on the line through $T_1$ other than $T_1R_i, T_1R, T_1G_5, T_1G$.
The final green point $G_6$ must be on $T_1R_j$.
If $G_6\in R_iG_2$, then clearly $R_i$ has ntype $[\ptype{1}{\ast}\ptype{1}{\ast}\ptype{1}{\ast}\ptype{4}{1}]$ and weight $6$.
If on the other hand $G_6\in R_iG_1$, then by Lemma~\ref{minimal} $R_iG_2$ is a $[\ptype{3}{1}]$-line; $R_i$ then has ntype $[\ptype{1}{4}\ptype{1}{b}\ptype{2}{c}\ptype{3}{1}]$, where $(b,c)=(4,2)$ or $(3,3)$, which gives $w(R_i)\in\{7/2,10/3\}$.
\end{proof}

\begin{claim}\label{b2}
Some $R_i$ has weight $6$.
\end{claim}
\begin{proof}[Proof of Claim~\ref{b2}]
Suppose not.
Since by Lemma~\ref{weight} the sum of the red weights must be $\binom{7}{2}$, we obtain by Claim~\ref{b1} that $\sum_{i=1}^4 w(R_i)\geq 12$.

Suppose $w(R_i)<3$ for some $i$, without loss for $i=4$.
Then $w(R_4)\leq 2$ and $\deg(R_4)\neq 4$ by Claim~\ref{b1}.
Since at most one red point has degree $4$ (Claim~\ref{claimtwo}), we have by Claim~\ref{b1}
\[ \sum_{i=1}^4 w(R_i) \leq 7/2 +3+3+2<12,\]
a contradiction.

Thus $w(R_i)\geq 3$ for all $i=1,2,3,4$.
By Claim~\ref{claimtwo} we have without loss one of the following two cases:
\begin{description}
\item[Case $\beta1$:] $w(R_i)=3$ for all $i=1,2,3,4$.
\item[Case $\beta2$:] $w(R_1)=w(R_2)=w(R_3)=3$, $w(R_4)\in\{10/3, 7/2\}$.
\end{description}

\subsection*{Case \boldmath$\beta1$\unboldmath}
Since $\sum_{i=1}^4 w(R_i)=12$ in this case, we must have $w(S_1)=w(S_2)=3/2$ by Claim~\ref{b1}.
If we now consider the subgeometry consisting of the green points, we find as in Case $\alpha$ that there are six $[\ptype{3}{\ast}]$-lines and three $[\ptype{2}{\ast}]$-lines, with some two of the $[\ptype{2}{\ast}]$-lines not intersecting in a green point, contradicting Lemma~\ref{sevenpointfls} as before.

\subsection*{Case \boldmath$\beta2$\unboldmath}
In this case $\sum_{i=1}^4 w(R_i)\in\{37/3, 25/2\}$, giving $w(S_1)+w(S_2)\in\{8/3,5/2\}$.
Without loss of generality we have one of 
\begin{itemize}
\item $w(S_1)=w(S_2)=4/3, w(R_4)=10/3$,
\item $w(S_1)=1, w(S_2)=3/2, w(R_4)=7/2$.
\end{itemize}
Considering the green subgeometry we obtain in both cases a contradiction as before.
\end{proof}

Thus without loss of generality $R_1$ has weight $6$ and ntype $[\ptype{1}{\ast}\ptype{1}{\ast}\ptype{1}{\ast}\ptype{4}{1}]$.
Denote the $[\ptype{4}{1}]$-line through $R_1$ by $\ell$.
Since $\deg(G)=5$, $GT_1$ must intersect $\ell$ in a green point, which we call $C$.
Let the other green point on $GT_1$ be $F$.
Choose a green point on $R_1T_1$ and call it $E$.
Then $GE$ intersects $\ell$ in a green point, say $D$.
Since $GE$ is then a $[\ptype{3}{1}]$-line, we must have $T_2\in GE$ (see Figure~\ref{fig17}).
\begin{figure}
\begin{center}
\begin{overpic}[scale=0.5]{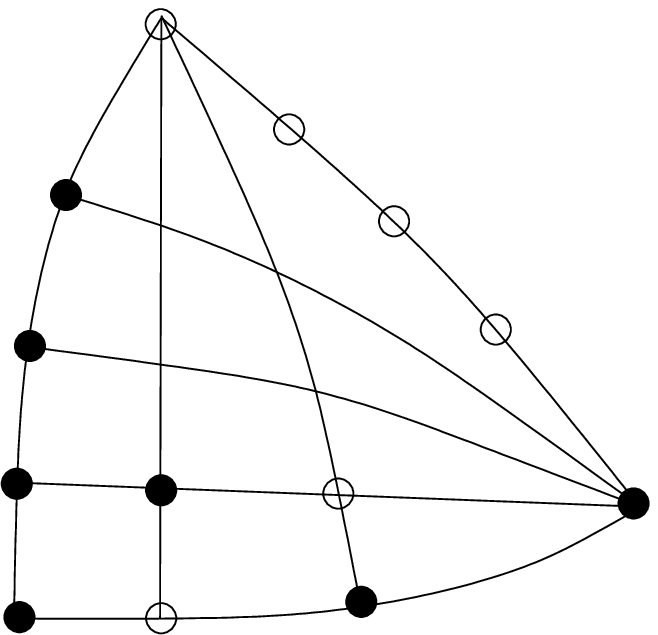}\small
\put(9,96){$R_1$}
\put(102,12){$G$}
\put(-3,66){$A$}
\put(-10,42){$B$}
\put(-12,21){$C$}
\put(-12,0){$D$}
\put(20,-10){$T_2$}
\put(51,-8){$E$}
\put(14,25.5){$F$}
\put(38,25){$T_1$}
\put(81,45){$R_4$}
\put(64,63){$R_3$}
\put(48,80){$R_2$}
\end{overpic}
\end{center}
\caption{The $(R_1,G)$-grid in Case $\beta$}\label{fig17}
\end{figure}
The only possible green point on $T_1T_2$ must be on $\ell$, say $A$.
Finally, we let $B$ be the fourth green point on $\ell$.
We thus have determined three full lines:
\[ CFGT_1, \quad DEGT_2, \quad R_1ABCD. \tag{2--4}\]
Without loss we may have $S_1\in AG$ and $S_2\in BG$.
There are now essentially two possible cases for the positions of $S_1$ and $S_2$.
\begin{description}
\item[Case $\beta3$] $R_1, S_1, S_2$ are not collinear.
\item[Case $\beta4$] $R_1, S_1, S_2$ are collinear.
\end{description}
We now show that $\beta3$ leads to a contradiction, while $\beta4$ leads to a unique MR geometry.

\subsection*{Case \boldmath$\beta3$\unboldmath}
Without loss of generality $S_1\in R_1T_1$ and $S_2\in R_1T_2$ (otherwise interchange $T_1\leftrightarrow T_2, C\leftrightarrow D, E\leftrightarrow F$).
It follows immediately that 
\[FR_1S_2T_2 \text{ and } ER_1S_1T_1\]
are full lines.
Since $T_1, T_2, A$ are collinear, $A\notin T_1S_2$.
Also, $B, C, E, F, G\notin T_1S_2$, and it follows that we must have $D\in T_1S_2$.
Similarly, since then $D\notin S_1S_2$, we obtain $C\in S_1S_2$.
Also, since then $C\notin S_1T_2$, we must have $B\in S_1T_2$.

We now consider the $(R_1, T_1)$-grid, recalling that the ntype of $T_1$ is known from Claim~\ref{b1} (Figure~\ref{fig19}).
\begin{figure}
\begin{center}
\begin{overpic}[scale=0.5]{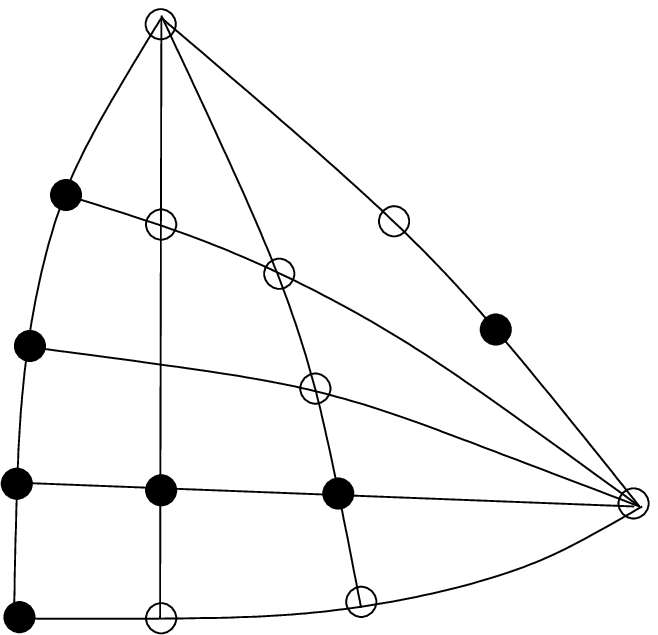}\small
\put(9,96){$R_1$}
\put(102,12){$T_1$}
\put(-3,66){$A$}
\put(-10,42){$B$}
\put(-12,21){$C$}
\put(-12,0){$D$}
\put(20,-10){$S_2$}
\put(51,-8){$R_4$}
\put(14,25.5){$F$}
\put(13,55){$T_2$}
\put(43,11){$G$}
\put(36,29){$R_3$}
\put(31,47){$R_2$}
\put(81,45){$E$}
\put(64,63){$S_1$}
\end{overpic}
\end{center}
\caption{The $(R_1,T_1)$-grid in Case $\beta3$}\label{fig19}
\end{figure}
We may assume without loss that $R_2\in AT_1, R_3\in BT_1, R_4\in DT_1$.
We now have determined three full lines: 
\[AR_2T_1T_2, \quad BR_3T_1, \quad DR_4S_2T_1.\]

Considering the $(R_1, T_2)$-grid (Figure~\ref{fig20}), we see that we must have $R_3\in CT_2$ (since $T_2\notin BR_3$) and then $R_4\in BT_2$.
Since $B\in S_1T_2$, we have determined the full lines
\[ BR_4S_1T_2 \text{ and } CR_3T_2.\]
\begin{figure}[b]
\begin{center}
\begin{overpic}[scale=0.5]{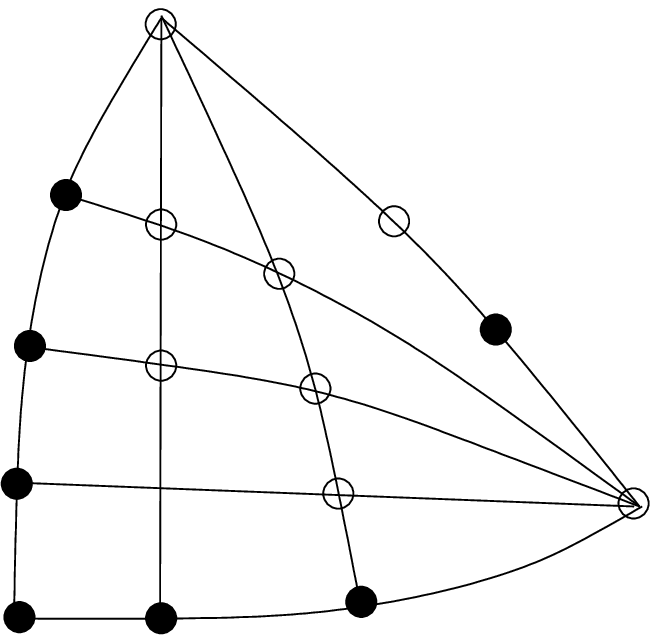}\small
\put(9,96){$R_1$}
\put(102,12){$T_2$}
\put(-3,66){$A$}
\put(-10,42){$B$}
\put(-12,21){$C$}
\put(-12,0){$D$}
\put(20,-10){$E$}
\put(51,-8){$G$}
\put(13,33){$S_1$}
\put(13,55){$T_1$}
\put(40,12){$R_3$}
\put(36,29){$R_4$}
\put(31,47){$R_2$}
\put(81,45){$F$}
\put(64,63){$S_2$}
\end{overpic}
\end{center}
\caption{The $(R_1,T_2)$-grid in Case $\beta3$}\label{fig20}
\end{figure}

By now we have the following full lines through $B$: $ABCDR_1$, $BGS_2$, $BS_1R_4T_2$, $BR_3T_1$.
Thus there must be exactly one more line through $B$, namely $BR_2$, which must also contain the remaining green points $E$ and $F$.
Thus \[BEFR_2\] is a full line.

Since $S_2\notin R_2A, R_2B, R_2G$, we obtain that $A,B,E,F,G\notin S_2R_2$.
Also, since $R_2\notin DS_2$, we have $D\notin S_2R_2$.
It follows that $C\in S_2R_2$.
We already have that $C, S_1, S_2$ are collinear, hence $C, S_1, S_2, R_2$ are collinear.
All other points are on the other lines through $R_2$ ($R_2A$, $R_2B$, $R_2G$).
Thus we have the full line \[CS_1S_2R_2,\] and $w(R_2)=3$.

We now have all the lines through $C$ except for $CR_4$, which must then be the full line \[CER_4.\]

We already have the following full lines through $R_3$: $GR_1R_2R_3R_4$, $BR_3T_1$, $CR_3T_2$.
Because of the full line $CS_1S_2R_2$ the points $S_1$, $S_2$, $R_3$ are not collinear.
Thus $R_3S_1$ and $R_3S_2$ are two more lines through $R_3$.
By Lemma~\ref{minimal} there must be a sixth line through $R_3$ as well.
Thus $R_3$ must have ntype $[\ptype{1}{4}\ptype{1}{2}\ptype{1}{2}\ptype{1}{2}\ptype{1}{2}\ptype{2}{1}]$ and $w(R_3)=1$.

Let $X_i$ be the green point on $R_3S_i$ for $i=1,2$, and let $X_3, X_4$ be the two green points on the $[\ptype{2}{1}]$-line through $R_3$.
Because of $AGS_1$ and $ER_1S_1T_1$, we have $X_1\neq A, E$, hence $X_1\in\{D,F\}$.
Because of $DR_4S_2T_1$ and $FR_1S_2T_2$, we have $X_2\neq D, F$, hence $X_2\in\{A,E\}$.

Because of $BEFR_2$, $\{X_3, X_4\}\neq\{E,F\}$, hence $(X_1, X_2)\neq (D, A)$.
Because of $DEGT_2$, $\{X_3, X_4\}\neq\{D,E\}$, hence $(X_1, X_2)\neq (F, A)$.
Because of $ABCDR_1$, $\{X_3, X_4\}\neq\{A,D\}$, hence $(X_1, X_2)\neq (F, E)$.
Therefore we must have $(X_1, X_2)=(D, E)$, hence $\{X_3, X_4\}=\{A, F\}$, and we have the full line
\[ AFR_3.\]

We already have the following full lines through $R_4$: $GR_1R_2R_3R_4$, $BR_4S_1T_2$, $DR_4S_2T_1$, $CER_4$.
Because of $ABCDR_1$, $R_4A$ is a fifth line.
Since $F\in AR_3$, $R_4F$ is a sixth line.
Thus $w(R_4)=1$.

If we now add up the weights of all the red points we obtain
\begin{eqnarray*}
21 &=& w(S_1)+w(S_2)+w(T_1)+w(T_2)\\
&& +w(R_1)+w(R_2)+w(R_3)+w(R_4)\\
&\leq& \frac32 + \frac32 +3+3+6+3+1+1 = 20,
\end{eqnarray*}
a contradiction.

\subsection*{Case \boldmath$\beta4$\unboldmath}
Without loss of generality $S_1, S_2\in R_1T_1$ (otherwise interchange $T_1\leftrightarrow T_2, C\leftrightarrow D, E\leftrightarrow F$).
It follows immediately that the following are full lines:
\[ FR_1T_2, \quad ER_1S_1S_2T_1, \quad AGS_1, \quad BGS_2.\tag{5--8}\]
Since $T_1, T_2, A$ are collinear, $A\notin T_2S_2$.
Also, $B, D, E, F, G\notin T_2S_2$, and it follows that we must have $C\in T_2S_2$.
Similarly, since then $C\notin T_2S_1$, we obtain $B\in T_2S_1$.

As in Case $\beta3$ we consider the $(R_1, T_1)$-grid (Figure~\ref{fig22}).
\begin{figure}
\begin{center}
\begin{overpic}[scale=0.5]{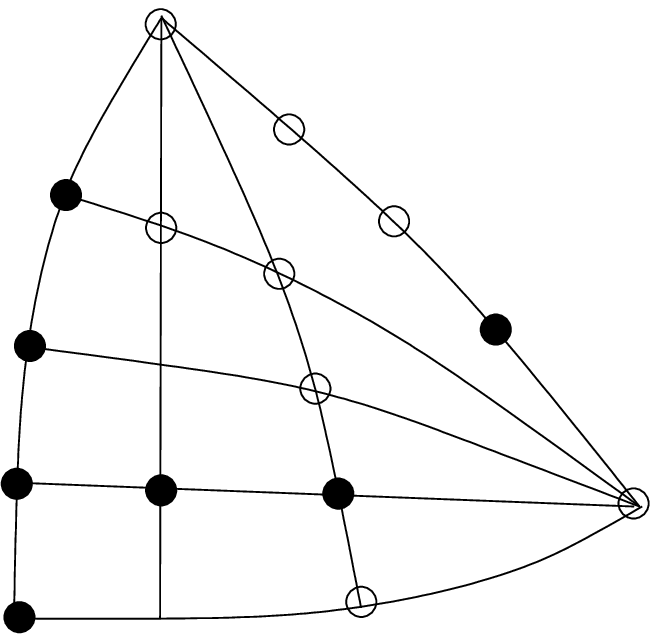}\small
\put(9,96){$R_1$}
\put(102,12){$T_1$}
\put(-3,66){$A$}
\put(-10,42){$B$}
\put(-12,21){$C$}
\put(-12,0){$D$}
\put(51,-8){$R_4$}
\put(14,25.5){$F$}
\put(81,45){$E$}
\put(64,63){$S_2$}
\put(48,80){$S_1$}
\put(13,55){$T_2$}
\put(43,11){$G$}
\put(36,29){$R_3$}
\put(31,47){$R_2$}
\end{overpic}
\end{center}
\caption{The $(R_1,T_1)$-grid in Case $\beta4$}\label{fig22}
\end{figure}
We may assume without loss that $R_2\in AT_1, R_3\in BT_1, R_4\in DT_1$.
We now have determined three full lines:
\[ AR_2T_1T_2, \quad BR_3T_1, \quad DR_4T_1.\tag{9--11}\]

From the $(R_1, T_2)$-grid (Figure~\ref{fig23}) we see that we must have $R_3\in CT_2$ (since $T_2\notin BR_3$), and then $R_4\in BT_2$, and also $S_2\in CT_2$ (since $B, S_1, T_2$ are collinear).
\begin{figure}[b]
\begin{center}
\begin{overpic}[scale=0.5]{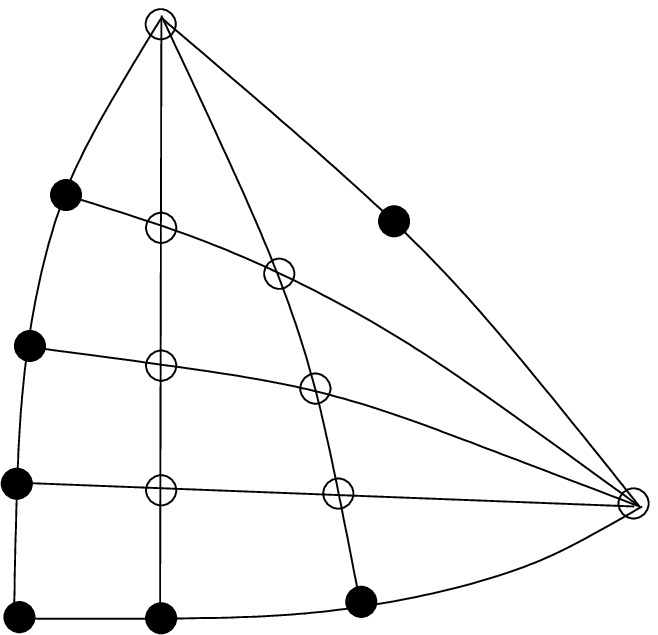}\small
\put(9,96){$R_1$}
\put(102,12){$T_2$}
\put(-3,66){$A$}
\put(-10,42){$B$}
\put(-12,21){$C$}
\put(-12,0){$D$}
\put(20,-10){$E$}
\put(51,-8){$G$}
\put(64,63){$F$}
\put(13,13){$S_2$}
\put(13,33){$S_1$}
\put(13,55){$T_1$}
\put(40,12){$R_3$}
\put(36,29){$R_4$}
\put(31,47){$R_2$}
\end{overpic}
\end{center}
\caption{The $(R_1,T_2)$-grid in Case $\beta4$}\label{fig23}
\end{figure}
We have determined the full lines
\[ BR_4S_1T_2, \quad CR_3S_2T_2. \tag{12--13}\]

We already have the following full lines through $B$: $ABCDR_1$, $BGS_2$, $BR_4S_1T_2$, $BR_3T_1$.
Thus there must be exactly one more line through $B$, namely $BR_2$, which must also contain the remaining green points $E$ and $F$.
Thus
\[ BEFR_2 \tag{14}\]
is a full line.

Now consider the $(R_1, E)$-grid (Figure~\ref{fig24}).
\begin{figure}
\begin{center}
\begin{overpic}[scale=0.5]{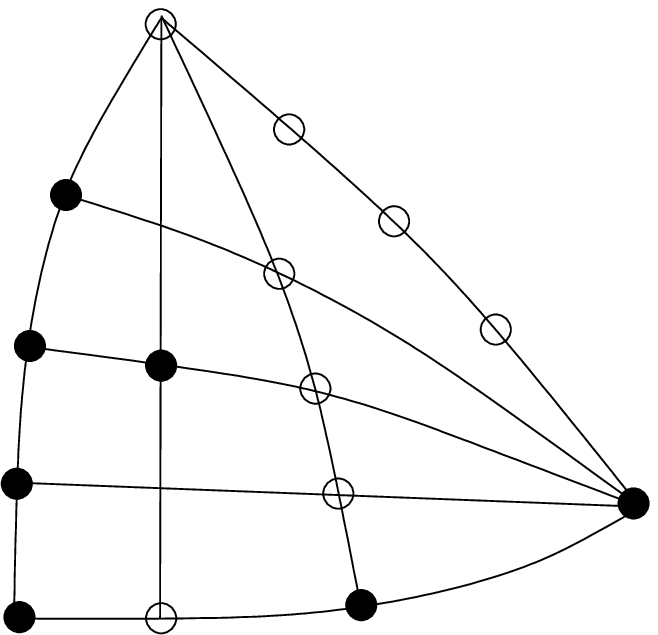}\small
\put(9,96){$R_1$}
\put(102,12){$E$}
\put(-3,66){$A$}
\put(-10,42){$B$}
\put(-12,21){$C$}
\put(-12,0){$D$}
\put(51,-8){$G$}
\put(13,33){$F$}
\put(81,45){$S_2$}
\put(64,63){$S_1$}
\put(48,80){$T_1$}
\put(20,-10){$T_2$}
\put(40,12){$R_4$}
\put(36,29){$R_2$}
\put(31,47){$R_3$}
\end{overpic}
\end{center}
\caption{The $(R_1,E)$-grid in Case $\beta4$}\label{fig24}
\end{figure}
There is a red point on each of $AE, BE, CE, DE$, and $4$ red points on $R_1E$.
Thus $\deg(E)=5$.
Since $E\notin CR_3, R_3\in AE$, and then $R_4\in CE$.
Thus we have the full lines
\[ AER_3, \quad CER_4.\tag{15--16}\]

Consider the $(R_1, R_2)$-grid (Figure~\ref{fig25}).
\begin{figure}[b]
\begin{center}
\begin{overpic}[scale=0.5]{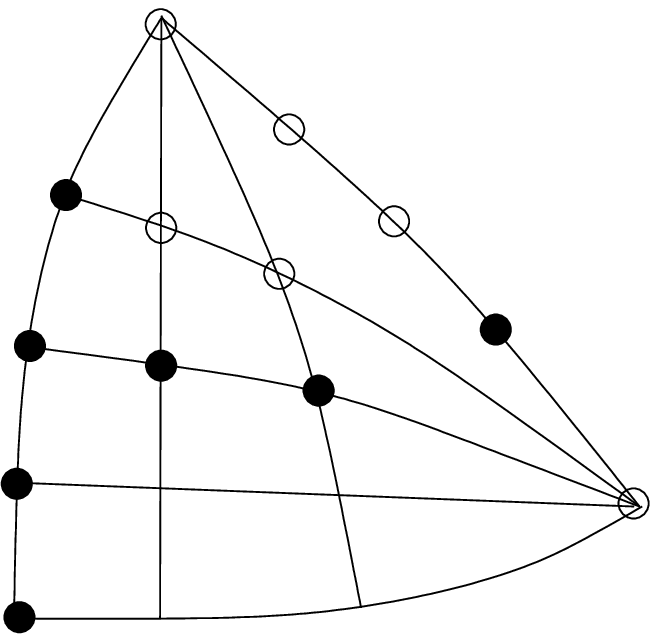}\small
\put(9,96){$R_1$}
\put(102,12){$R_2$}
\put(-3,66){$A$}
\put(-10,42){$B$}
\put(-12,21){$C$}
\put(-12,0){$D$}
\put(13,33){$E$}
\put(81,45){$G$}
\put(64,63){$R_4$}
\put(48,80){$R_3$}
\put(36,29){$F$}
\put(31,47){$T_2$}
\put(13,55){$T_1$}
\end{overpic}
\end{center}
\caption{The $(R_1,R_2)$-grid in Case $\beta4$}\label{fig25}
\end{figure}
Considering the green points we see that $\deg(R_2)=5$ and $w(R_2)=3$.
Since $R_2\notin CS_2$, $S_2\notin CR_2$, hence $S_2\in DR_2$ and $S_1\in CR_2$.
Thus we have the full lines
\[ CR_2S_1, \quad DR_2S_2.\tag{17--18}\]

Consider the $(R_1, R_3)$-grid (Figure~\ref{fig26}).
\begin{figure}
\begin{center}
\begin{overpic}[scale=0.5]{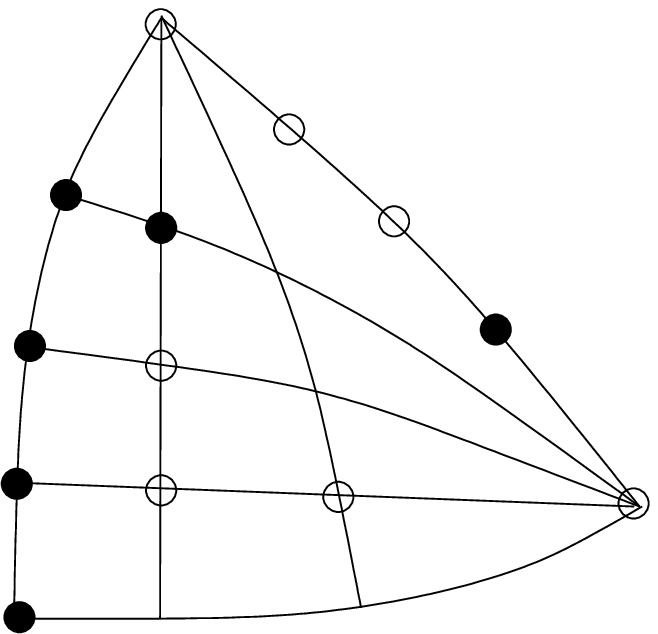}\small
\put(9,96){$R_1$}
\put(102,12){$R_3$}
\put(-3,66){$A$}
\put(-10,42){$B$}
\put(-12,21){$C$}
\put(-12,0){$D$}
\put(13,33){$T_1$}
\put(81,45){$G$}
\put(64,63){$R_4$}
\put(48,80){$R_2$}
\put(13,13){$S_2$}
\put(40,12){$T_2$}
\put(13,55){$E$}
\end{overpic}
\end{center}
\caption{The $(R_1,R_3)$-grid in Case $\beta4$}\label{fig26}
\end{figure}
If $S_1\in DR_3$, then either $F\in DR_3$ and $w(R_3)=3/2$, or $F\notin DR_3$ and $w(R_3)=1$.
If on the other hand $S_1\notin DR_3$, then $F\in S_1R_3$ and $w(R_3)=1$.
Summarizing: $w(R_3)\leq 3/2$, with equality iff $DFR_3S_1$ is a line.

Consider the $(R_1, R_4)$-grid (Figure~\ref{fig27}).
\begin{figure}[b]
\begin{center}
\begin{overpic}[scale=0.5]{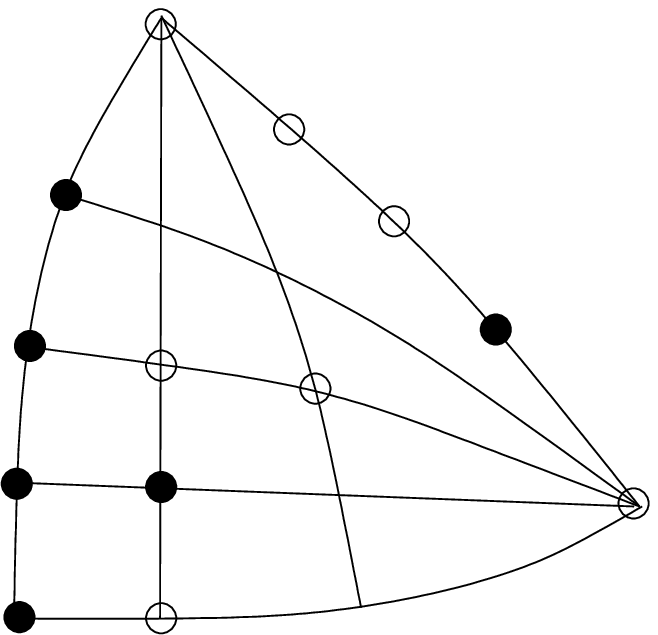}\small
\put(9,96){$R_1$}
\put(102,12){$R_4$}
\put(-3,66){$A$}
\put(-10,42){$B$}
\put(-12,21){$C$}
\put(-12,0){$D$}
\put(13,33){$S_1$}
\put(81,45){$G$}
\put(64,63){$R_3$}
\put(48,80){$R_2$}
\put(13,13){$E$}
\put(20,-10){$T_1$}
\put(36,29){$T_2$}
\end{overpic}
\end{center}
\caption{The $(R_1,R_4)$-grid in Case $\beta4$}\label{fig27}
\end{figure}
If $S_2\in AR_4$, then either $F\in AR_4$ and $w(R_4)=3/2$, or $F\notin AR_4$ and $w(R_4)=1$.
If on the other hand $S_2\notin AR_4$, then $F\in S_2R_4$ and $w(R_4)=1$.
Summarizing: $w(R_4)\leq 3/2$, with equality iff $AFR_4S_2$ is a line.

We now add up the weights of the red points:
\begin{eqnarray*}
21 &=& w(S_1)+w(S_2)+w(T_1)+w(T_2)\\
&& +w(R_1)+w(R_2)+w(R_3)+w(R_4)\\
&\leq& \frac32 + \frac32 +3+3+6+3 + \frac32 + \frac32 = 21.
\end{eqnarray*}
It follows that $w(R_3)=w(R_4)=3/2$, giving the lines
\[ DFR_3S_1, \quad AFR_4S_2.\tag{19--20}\]
We have now determined all the lines (1--20).
We have found a unique MR geometry, which is isomorphic to $\MR_{15}^5$ by the following isomorphism:
\begin{center}
\begin{tabular}{c|c|c|c|c|c|c|c|c|c|c|c|c|c|c}
$A$ & \!$B$ & \!$C$ & \!$D$ & \!$E$ & \!$F$ & \!$G$ & \!$R_1$ & \!$R_2$ & \!$R_3$ & \!$R_4$ & \!$S_1$ & \!$S_2$ & \!$T_1$ & \!$T_2$ \\ \hline
$5_2$ & \!$4_0$ & \!$5_1$ & \!$3_0$ & \!$0_0$ & \!$1_2$ & \!$1_0$ & \!$6_0$ & \!$5_0$ & \!$2_2$ & \!$2_1$ & \!$4_2$ & \!$2_0$ & \!$4_1$ & \!$6_1$
\end{tabular}
\end{center}
We have finished the proof of Proposition~\ref{prop6} and thus of Theorem~\ref{mainthm}.
\end{proof}

\section*{Appendix: Finding all one-point extensions of the MR geometries on \boldmath$14$\unboldmath\ points}
A \emph{partial parallel class} (ppc) of a finite linear space $S$ is a partition of the points of $S$ such that each set in the partition is either a line of $S$ or a singleton.
An \emph{extendable} ppc (eppc) of an MR geometry is a ppc in which all singletons (if any) are of the same colour.
Clearly each one-point extension of an f.l.s.\ $S$ determines a ppc:
If $S\cup\{\infty\}$ is an extension of $S$ then
\[ \{\ell\setminus\{\infty\} : \infty\in\ell\} \]
is a ppc of $S$.
Furthermore, for each MR geometry, this ppc is an eppc.
Conversely, each eppc of an MR geometry determines at least one one-point extension:
\begin{lemma}\label{ppclemma}
Let $(S,\chi)$ be an MR geometry.
\begin{itemize}
\item For each ppc of $(S,\chi)$ not containing singletons there exist exactly two one-point extensions of $(S,\chi)$.
\item For each ppc of $(S,\chi)$ containing at least one singleton, and with all singletons of the same colour, there exits exactly one one-point extension of $(S,\chi)$.
\end{itemize}
\end{lemma}
\begin{proof}
If $\{B_i: i\in I\}$ is a ppc then the extension we are looking for is $S\cup\{\infty\}$ with the original lines of $S$ together with all $B_i\cup\{\infty\}$.
In the first case we are free to give $\infty$ any colour, while in the second case the colour of $\infty$ is forced.
\end{proof}

\subsection*{Extending \boldmath$\MR_{14}^1$\unboldmath}
By Lemma~\ref{ppclemma} it is sufficient to enumerate all eppc's of $\MR_{14}^1$.

\subsubsection*{\boldmath Case I: The eppc contains a $4$-line of $\MR_{14}^1$\unboldmath}
By symmetry we may take the line to be $0_01_03_06_1$.
The only lines parallel to $0_01_03_06_1$ are $0_11_13_16_0$, $5_05_1$, $2_02_1$, and $4_04_1$.
These five lines already partition $\MR_{14}^1$.
If one of these lines, say $\ell$, does not belong to the eppc then the eppc must contain singletons, one corresponding to each point of $\ell$.
Then there are singletons of different colours, a contradiction.
Thus the only eppc is
\[ \{0_01_03_06_1, 0_11_13_16_0, 5_05_1, 2_02_1, 4_04_1\}.\]
The point at infinity on these lines correspond to $6_2$ in $\PG(2,4)$.
Since the eppc has no singletons, $6_2$ can be assigned any colour.
Thus there are two extensions that are clearly isomorphic.
If we colour $6_2$ green we obtain $\MR_{15}^1$.

\subsubsection*{\boldmath Case II: There are no $4$-lines in the eppc\unboldmath}
The seven $2$-lines of $\MR_{14}^1$ already form a partition.
If one of them does not belong to the eppc, then as before, there are singletons of different colours, a contradiction.

Thus there are two one-point extensions, depending on the colour of $\infty$ adjoined to each of the $2$-lines.
These two extensions are clearly isomorphic.
If we colour $\infty$ green we obtain $\MR_{15}^2$.

\subsection*{Extending \boldmath$\MR_{14}^2$\unboldmath}
Since it is not clear from the given representation of $\MR_{14}^2$ what its symmetries are, we use another representation:
\begin{center}
\begin{tabular}{ll}
Green points: & $G, G_1, G_2, G_3, H_1, H_2, H_3$\\
Red points: & $R, R_1, R_2, R_3,, S_1, S_2, S_3$ \smallskip\\
Lines: & $GR_1R_2R_3R$, $RH_1H_2H_3$\\
& $GG_iH_iS_i$, $G_iG_{i+1}H_{i-1}R_{i+1}$\\
& $H_{i-1}R_{i-1}S_iS_{i+1}$, $RS_iG_{i-1}$\\
& $G_{i+1}S_iR_i$, $H_iR_{i+1}$,
\end{tabular}
\end{center}
with $i=1,2,3$, and indices modulo $3$.

The mappings in Figure~\ref{fig28} can easily be checked to be isomorphisms (between linear spaces).
\begin{figure}
\begin{center}
\begin{tabular}{c|c|c|c|c|c|c|c|c|c|c|c|c|c}
$G$ & \!$R$ & \!$G_1$ & \!$G_2$ & \!$G_3$ & \!$S_1$ & \!$S_2$ & \!$S_3$ & \!$H_1$ & \!$H_2$ & \!$H_3$ & \!$R_1$ & \!$R_2$ & \!$R_3$ \\ \hline\hline
$6_0$ & \!$6_1$ & \!$4_1$ & \!$2_1$ & \!$5_1$ & \!$2_0$ & \!$5_0$ & \!$4_0$ & \!$0_0$ & \!$1_0$ & \!$3_0$ & \!$1_2$ & \!$3_2$ & \!$0_2$ \\ \hline
$6_0$ & \!$6_1$ & \!$2_0$ & \!$4_0$ & \!$5_0$ & \!$4_1$ & \!$5_1$ & \!$2_1$ & \!$0_0$ & \!$3_0$ & \!$1_0$ & \!$0_2$ & \!$3_2$ & \!$1_2$ \\ \hline
$6_0$ & \!$6_1$ & \!$5_0$ & \!$2_0$ & \!$4_0$ & \!$2_1$ & \!$4_1$ & \!$5_1$ & \!$1_0$ & \!$0_0$ & \!$3_0$ & \!$1_2$ & \!$0_2$ & \!$3_2$ \\ \hline
$6_0$ & \!$6_1$ & \!$2_1$ & \!$5_1$ & \!$4_1$ & \!$5_0$ & \!$4_0$ & \!$2_0$ & \!$1_0$ & \!$3_0$ & \!$0_0$ & \!$3_2$ & \!$0_2$ & \!$1_2$ \\ \hline
$6_0$ & \!$6_1$ & \!$5_1$ & \!$4_1$ & \!$2_1$ & \!$4_0$ & \!$2_0$ & \!$5_0$ & \!$3_0$ & \!$0_0$ & \!$1_0$ & \!$0_2$ & \!$1_2$ & \!$3_2$ \\ \hline
$6_0$ & \!$6_1$ & \!$4_0$ & \!$5_0$ & \!$2_0$ & \!$5_1$ & \!$2_1$ & \!$4_1$ & \!$3_0$ & \!$1_0$ & \!$0_0$ & \!$3_2$ & \!$1_2$ & \!$0_2$
\end{tabular}
\end{center}
\caption{All isomorphisms between the two representations of $\MR_{14}^2$}\label{fig28}
\end{figure}
We now show that these are the only isomorphisms.

The only point occurring in only one $4$-line is $6_1$ in the first representation, and $R$ in the second representation.
Thus we must have $R\mapsto 6_1$ and $\{H_1,H_2,H_3\}\mapsto\{0_0,1_0,3_0\}$.

The only $2$-lines are $1_23_0$, $3_20_0$, $0_21_0$, and $H_1R_2$, $H_2R_3$, $H_3R_1$.
Thus we must have $\{R_1,R_2,R_3\}\mapsto\{0_2,1_2,3_2\}$.
Furthermore, the restriction to $\{H_1,H_2,H_3\}$ determines the restriction to $\{R_1,R_2,R_3\}$.

The only $5$-line is $0_21_23_26_06_1$, $GR_1R_2R_3R$.
Thus we must have $G\mapsto6_0$.

At this stage $G$ and $R$ are fixed.
If we choose any of the $3!$ mappings $\{H_1,H_2,H_3\}\mapsto\{0_0,1_0,3_0\}$ then, as mentioned above, the images of $R_1, R_2, R_3$ are determined.
In each case, by considering the three lines $H_iR_i$, we determine the images of the pairs $\{S_2,S_3\}$, $\{S_3,S_1\}$, $\{S_1,S_2\}$, from which the image of each $S_i$ is determined.
Then by considering the lines $RS_i$ we determine the images of $G_1, G_2, G_3$, and the whole isomorphism is fixed.
It then only remains to check that the remaining lines $GG_iH_iS_i$, $G_iG_{i+1}H_{i-1}R_{i+1}$, $G_{i+1}S_iR_i$ are mapped properly.
In this way all six isomorphisms are obtained.
(Note that the even automorphisms preserve the colours, i.e., they are also isomorphisms of MR geometries.)

Now we can write down the automorphisms of the second representation:
\begin{center}
\begin{tabular}{l}
identity\\
$(G_1S_1)(G_2S_3)(G_3S_2)(H_3H_2)(R_3R_1)$\\
$(G_1S_2)(G_2S_1)(G_3S_3)(H_2H_1)(R_2R_3)$\\
$(G_1S_3)(G_2S_2)(G_3S_1)(H_1H_3)(R_1R_2)$\\
$(G_1G_2G_3)(S_1S_2S_3)(H_1H_2H_3)(R_1R_2R_3)$\\
$(G_3G_2G_1)(S_3S_2S_1)(H_3H_2H_1)(R_3R_2R_1)$
\end{tabular}
\end{center}
The action of the automorphism group on the lines has seven orbits:
\begin{center}
\begin{tabular}{ll}
(i) & $GR_1R_2R_3R$\\
(ii) & $RH_1H_2H_3$\\
(iii) & $GG_iH_iS_i$\\
(iv) & $G_iG_{i+1}H_{i-1}R_{i+1}$, $H_{i-1}R_{i-1}S_iS_{i+1}$\\
(v) & $RS_iG_{i-1}$\\
(vi) & $G_{i+1}S_iR_i$\\
(vii) & $H_iR_{i+1}$
\end{tabular}
\end{center}

We now determine all eppc's of $\MR_{14}^2$.
We consider seven cases, depending on whether the eppc contains a line from one of the seven orbits.
(We could also have considered $20$ cases, one for each line, in which case it wouldn't have been necessary to calculate the automorphism group.
However, such an approach would have been more tedious and less insightful.)

\subsubsection*{Case (i)}
Since $GR_1R_2R_3R$ intersects all other lines, all points not on this line must be singletons in a ppc.
Thus there is no eppc containing $GR_1R_2R_3R$.

\subsubsection*{Case (ii)}
Suppose $RH_1H_2H_3$ is an eppc.
The only parallel lines are $G_{i+1}S_iR_i$.
All of them must be in the eppc.
The only remaining point is $G$, which must be a singleton in the eppc.
Thus there is one extension, with $\infty$ coloured red.
This MR geometry is isomorphic to $\MR_{15}^1$:
\begin{center}
\begin{tabular}{c|c|c|c|c|c|c|c|c|c|c|c|c|c|c}
$G$ & \!\!$R$ & \!\!$G_1$ & \!\!$G_2$ & \!\!$G_3$ & \!\!$S_1$ & \!\!$S_2$ & \!\!$S_3$ & \!\!$H_1$ & \!\!$H_2$ & \!\!$H_3$ & \!\!$R_1$ & \!\!$R_2$ & \!\!$R_3$ & \!\!$\infty$ \\ \hline
$6_1$ & \!$6_2$ & \!$5_1$ & \!$4_1$ & \!$2_1$ & \!$2_0$ & \!$5_0$ & \!$4_0$ & \!$1_1$ & \!$3_1$ & \!$0_1$ & \!$0_0$ & \!$1_0$ & \!$3_0$ & \!$6_0$
\end{tabular}
\end{center}

\subsubsection*{Case (iii)}
Assume without loss of generality that $RH_1H_2H_3$ is not in the eppc, and $GG_1H_1S_1$ is.
Lines parallel to $GG_1H_1S_1$ are $RS_3G_2$, $G_3S_2R_2$, $H_2R_3$, $H_3R_1$.
They already partition $\MR_{14}^2$.
All of them must be in the eppc.
The point $\infty$ may have either colour.
In the first two isomorphisms in Figure~\ref{fig28} the point $\infty$ corresponds to $4_2$.
We have obtained $\MR_{15}^{3\mathsf{g}}$ and $\MR_{15}^{3\mathsf{r}}$.

\subsubsection*{Case (iv)}
Assume without loss of generality that Cases (i), (ii) and (iii) do not occur, and that $G_1G_2H_3R_2$ is in an eppc.
Lines parallel to $G_1G_2H_3R_2$ are $H_1R_1S_2S_3$, $RS_1G_3$, $H_2R_3$.
All these lines must be in the eppc, leaving $G$ as a singleton.
Thus there is a single extension, with $\infty$, which has to be red, corresponding to $3_1$ in the second and fourth isomorphisms in Figure~\ref{fig28}.
We have obtained $\MR_{15}^4$.

\subsubsection*{Cases (v), (vi), and (vii)}
In each of the remaining cases we obtain that a ppc must contain singletons of different colours, giving no extensions.

\bigskip
We have shown that  $\MR_{15}^1$, $\MR_{15}^{3\mathsf{g}}$, $\MR_{15}^{3\mathsf{r}}$, and $\MR_{15}^4$ are the only extensions of $\MR_{14}^2$.\qed

\section*{Acknowledgements}
We thank the referee for valuable advice.

\end{document}